   \def\obs#1{{\bf (*** #1 ***)} }

 \def\obs#1{}     
\documentclass[twoside,letterpaper,draft,11pt]{amsart}

\usepackage{amsmath,amsthm,latexsym,xspace,amscd,amssymb,xypic}               
\usepackage[dvipsnames]{xcolor}

\title[ A  seven-term exact sequence]{Partial generalized crossed products and a  seven-term exact sequence}

\author[M.\ Dokuchaev]{M. Dokuchaev}
\address{Instituto de
Matem\'atica e Estat\'\i stica,
Universidade de S\~ao Paulo, Rua do Mat\~ao, 1010,
05508-090 S\~ao Paulo, SP, Brasil.}
\email{dokucha@ime.usp.br}

\author[A.\ Paques]{A. Paques}
\address{Instituto de
Matem\'atica  Universidade  Federal do Rio Grande do Sul,  Avenida Bento Gon\c calves 9500,
91509-900 Porto Alegre, RS, Brasil.}
\email{paques@mat.ufrgs.br}

\author[H.\ Pinedo ]{H. Pinedo }
\address{Escuela de
Matem\'aticas, Universidad Industrial de Santander, Cra 27 calle 9, Edificio Camilo Torres, Bucaramanga, Colombia.}
\email{hpinedot@uis.edu.co}

\author[I.\ Rocha]{I. Rocha}
\address{Unidade Acad\^emica de Matem\'atica, Universidade Federal de Campina Grande, Avenida Apr\'{\i}gio Veloso, 882, 58429-900, Campina Grande, PB, Brasil.}
\email{itailma@mat.ufcg.edu.br}

\thanks{{This work was partially supported by FAPESP of Brazil. The first named author was also partially supported by CNPq of Brazil.}\\{\bf  Mathematics Subject Classification}:
Primary 13B05; Secondary   13A50; 16H05; 16K50; 16S35; 16W22; 20M18.\\
{\bf Key words and phrases:} {Partial action, partial representation,  Galois extension, Galois cohomology, Azumaya algebra, Brauer group, Picard group.}}

\newtheorem{teo}{Theorem}[section]
\newtheorem{defi}[teo]{Definition}
\newtheorem{lema}[teo]{Lemma}
\newtheorem{cor}[teo]{Corollary}
\newtheorem{prop}[teo]{Proposition}

\newtheorem{exe}[teo]{Example}
\newtheorem{remark}[teo]{Remark}

\newcommand{\La}{\Lambda}

\newcommand{\0}{\theta}

\newcommand{\de}{\delta}

\newcommand{\cua}{\hfill${\blacksquare}$}
\newcommand{\omgh}{\omega_{g,h}}
\newcommand{\af}{\alpha}
\newcommand{\afg}{\alpha_g}
\newcommand{\afh}{\alpha_h}

\newcommand{\afgh}{\alpha_{gh}}

\newcommand{\Pics}{\mathbf{PicS}}
\newcommand{\lb}{\lambda}

\newcommand{\te}{\theta}

\newcommand{\p}{{\bf Proof. }}

\newcommand{\B}{{\mathcal B}}

\newcommand{\mfp}{{\mathfrak p}}

\newcommand{\U}{{\mathcal U}}
\newcommand{\df}{\displaystyle\frac}
\newcommand{\ot}{{\otimes}}

\newcommand{\di}{{\diamond}}
\newcommand{\Pic}{\mathbf{Pic}}
\newcommand{\m}{{}^{-1}}
\newcommand{\mt}{\mapsto}
\newcommand{\om}{\omega}
\def\f{\varphi}

\def\ndv{\ {\mid \kern -0.7 em {\scriptstyle \not}} \ \ }

\def\nd{\ {\mid \kern -0.4 em {\scriptstyle \not}} \ \ }

\newcommand{\N}{{\mathbb N}}

\newcommand{\vu}{\vspace{.1cm}}
\newcommand{\vd}{\vspace{.2cm}}
\newcommand{\vt}{\vspace{.3cm}}

\begin{document}

\date{\today}

\begin{abstract}  For a partial Galois extension of commutative rings we give a seven term exact sequence wich
generalize the  Chase-Harrison-Rosenberg  sequence.
\end{abstract}
\maketitle
\begin{section}{Introduction}

This is a continuation of the paper \cite{DPP1}, in which a sequence of six homomorphisms was constructed in order to give a partial Galois theoretic analogue of the Chase-Harrison-Rosenberg seven term exact sequence \cite{CHR}. In this article we prove that the sequence from  \cite{DPP1} is exact. In some sense this topic is  a natural sequel to
\cite{DFP}, where a Galois theory for commutative rings was  developed in the context of partial actions of groups on rings, including some facts on non-commutative Galois extensions. A systematic algebraic consideration of partial actions was initially influenced by a successful use of partial actions and partial representations  in the study of $C^*$-algebras generated by partial isometries, and nowadays it is being stimulated by diverse algebraic, topological and  $C^*$-algebraic advances  on the subject (see the book by R. Exel \cite{E6} and  the surveys \cite{Ba}, \cite{D2}, \cite{D3}, \cite{F2}, \cite{Pi3}).\\

Prior to the sequence, the paper by S.U. Chase,  D.K. Harrison, A. Rosenberg  \cite{CHR}  contains,
 in particular,    several equivalent definitions of the Auslander-Goldman concept of a Galois extension of commutative rings \cite{AG} and a fundamental theorem. This part of \cite{CHR} was extended in  \cite{DFP} to partial actions, naturally raising the problem of the generalization of the seven term exact sequence as well. Given  a (global) Galois extension $R^G\subseteq R$ of commutative rings relative to a finite group $G,$ the Chase-Harrison-Rosenberg  sequence involves Picard groups,  Galois cohomology groups and the   Auslander-Goldman notion of the  Brauer group   \cite{AG}, which in this case is formed by  the equivalence classes of Azumaya (i.e. central separable) $R^G$-algebras split  by $R.$ In symbols,  the sequence is of the form
\begin{align*}
&0\to H^1(G, \U(R)) {\to} {\bf Pic}(R^G) {\to}  {\bf Pic}(R)^G {\to }
H^2(G, \U(R)  ) {\to} B(R/R^G) {\to} H^1(G,{\bf Pic}(R)){\to} \\ &H^3(G, \U(R)).
\end{align*}
Thus for our generalization of  the sequence a cohomology theory  based on partial group actions is needed. The latter was introduced in \cite{DK}, inspired by R. Exel's notion of a twisted partial group action on a $C^*$-algebra \cite{E0}, adapted to rings in \cite{DES}.   The starting point of  \cite{DK} is the replacement of global $G$-actions on abelian groups ($G$-modules) by unital partial actions of $G$ on commutative semigroups (partial $G$-modules), in par\-ti\-cu\-lar, on commutative rings.   The partial  group cohomology  found applications to partial projective group representations   \cite{DK}, \cite{DoSa}
 and  to the study of ideals of (global)  reduced $C^*$-crossed products   \cite{KennedySchafhauser}. It also  motivated  the treatment of partial cohomology from the point of view of  Hopf algebras  \cite{BaMoTe}.  A more general multiplier valued version of partial cohomology was used to classify  extensions of semilattices of abelian groups  by groups   \cite{DKh2}, \cite {DKh3}.\\

The action of $G$ on $R,$ as above, induces an action of $G$ on the Picad group ${\bf Pic}(R),$ and the third term of the   Chase-Harrison-Rosenberg sequence is the subgroup of the fixed points  ${\bf Pic}(R)^G.$ It turns out that a partial action $\alpha $  of $G$ on $R$ gives rise to a partial action $\alpha ^*$ of $G$ on a more general object, namely the commutative inverse semigroup
${\bf PicS}(R)$ of the finitely generated projective $R$-modules of rank $\leq 1$  (see Section~\ref{f_1f_2f_3}).  The semigroup ${\bf PicS}(R)$ contains  ${\bf Pic}(R)$ as a subgroup, and  the appropriate replacement of   ${\bf Pic}(R)^G$ is the semigroup of those elements of  ${\bf Pic}(R),$ which are  $\alpha ^*$-invariant in  ${\bf PicS}(R).$ All (global) cohomology groups in the above sequence are substituted now by their partial theoretic analogues, observing that ${\bf PicS}(R)$ appears also in the sixth term as the coefficients of certain  $1$-cohomologies. Our main result is the following exact sequence:\\
\begin{align*}
&0\to H^1(G,\af, R)\stackrel{\f_1}{\to} {\bf Pic}(R^\af)\stackrel{\f_2}{\to }{\bf PicS}(R)^{\af^*}\cap {\bf Pic}(R) \stackrel {\f_3}{\to }H^2(G,\af, R) \stackrel{\f_4}{\to} B(R/R^\af)\stackrel{\f_5}{\to}\\& H^1(G,\af^*,{\bf PicS}(R))\stackrel{\f_6}{\to}   H^3 (G,\af, R).
\end{align*}

The seven term sequence in   \cite{CHR} was derived from the Amitsur cohomology  seven term exact sequence, obtained by S.U. Chase and A. Rosenberg in  \cite{CR} using a spectral sequence of A. Grothendieck. A constructive proof of the Chase-Harrison-Rosenberg sequence  was offered by T. Kanzaki in \cite{K}, an important ingredient being the notion of a generalized crossed product related to a  representation of $G$ into an appropriate  Picard group. Our approach is   inspired by the treatment of the first six terms of the sequence in the  book of F. DeMeyer and  E. Ingraham  \cite{DI} and the Kanzaki's paper \cite{K}. We use an extended version of generalized crossed products, wich in our case are related to a partial representation of $G$ into a certain Picard semigroup (see Section~\ref{pargencrosprod}).\\

In this paper by a ring we mean an associative ring with identity element. For any ring $R,$ an $R$-module is a left unital  $R$-module.
If $R$ is commutative, unless otherwise stated,  we will consider an $R$-module $M$ as a central $R$-$R$-bimodule.
 We say that $M$ is a {\it  f.g.p. $R$-module} if $M$ is a (left) projective and finitely generated $R$-module. A
 {\it  faithfully projective $R$-module} is a faithful f.g.p. $R$-module. For a monoid (or a ring)
 $T,$ the group  of units of $T$ is denoted by $\U(T).$\\

Throughout this work, $R$ will denote a commutative ring and unadorned $\ot,$ unless otherwise established, will mean $\ot_R.$

\end{section}

\begin{section}{Some preliminaries}

In this section we give some definitions and results which will be used  in  the paper.
We start by recalling three basic facts about projective modules.
\begin{lema}  The following assertions hold.
\begin{itemize}
\item \label{proje}\cite[{ Lemma 2.1}]{DPP1} Let  $M$ and $N$ be R-modules such that $M$ and  $M\ot N$ are f.g.p. R-modules, and
$M_\mfp\ne 0$ for all prime ideals $\mfp$ of $R.$ Then, N is also a  f.g.p. R-module.
\item
\cite[Lemma I.3.2 b)]{DI} If $M$ is a faithful f.g.p. R-module, then there exists an\\ $R$-$R$-bimodule isomorphism
$M^*\ot_{  {  {\rm End}_{R}  }  (M)} M\cong R.$ Consequently, if $N$ and $N'$ are $R$-modules such that
$M\ot N\cong M\ot N'$ as $R$-modules, then $N\cong N'$ as $R$-modules.

\item \label{ht}\cite[Hom-Tensor Relation I.2.4]{DI}  Let  $A$ and $B$ be $R$-algebras, $M$ be a f.g.p. $A$-module and
$N$ be a f.g.p. $B$-module. Then for any $A$-module $M'$ and any $B$-module $N'$, the map
 $$\psi\colon {\rm Hom}_A(M,M')\ot {\rm Hom}_B(N,N')\to {\rm Hom}_{(A\ot B)}(M\ot N,M'\ot N' ),$$
 induced by $(f\ot  g)(m\ot n)=f(m)\ot g(n),$ for all $m\in M,\, n\in N$, is an $R$-module isomorphism. If $M=M'$ and $N=N',$
 then $\psi$ is an $R$-algebra isomorphism.
\end{itemize}
\end{lema}

 We proceed with partial Galois cohomology.
Let $G$ be a group.
 We recall from \cite{DES} that a {\it unital twisted partial action} of  $G$ on  $R$ is a triple $$\af=(\{D_g\}_{g\in G}, \{\afg\}_{g\in G}, \{\omgh\}_{(g,h)\in G\times G}),$$
such that for every $g\in G,$ $D_g$ is an ideal of $R$ generated by a non-necessarily non-zero idempotent $1_g,$ $\afg\colon D_{g\m}\to D_g$ is
an isomorphism of rings,  for each pair $(g,h)\in G\times G,$ $\omgh$ is in  $\U(D_gD_{gh}),$ and the following statements are satisfied
for all $g,h,l\in G:$

\begin{itemize}

\vu

\item [(i)] $D_1=R$ and $\af_1$ is the identity map of $R,$

\vu

\item [(ii)] $\afg(D_{g\m}D_h)=D_gD_{gh},$

\vu

\item [(iii)] $\afg \circ \afh(r)=\omgh\afgh(r)\omgh^{-1}$ for any $r\in D_{h\m} D_{(gh)\m},$

\vu

\item [(iv)] $\om_{1,g}=\om_{g,1}=1_g$  and

\vu

\item [(v)] $\afg(1_{g\m}\om_{h,l})\om_{g,hl}=\omgh\om_{gh,l}.$
\end{itemize}

\vu

If $R$ is a multiplicative monoid,  one obtains  from the above definition the concept of a unital twisted partial action of a group on a monoid.

\vd

Notice that (v) implies
\begin{equation}\label{afom}\afg(\om_{g\m,g})=\om_{g,g\m},\,\, \,\,\text{for any}\,\,g\in G.\end{equation}

\vd

Let $\af$ be a unital twisted partial action of $G$ on $R.$ The {\it partial crossed product} $R\star_{\af,\om} G$ (see \cite[Definition 2.2]{DES}) is the abelian group  $\bigoplus\limits_{g\in G}D_g\delta_g,$ where the $\delta_g$' s are place holder symbols, and the multiplication
in $R\star_{\af,\om} G$ is induced  by \begin{equation*}\label{product}(r_g\delta_g)(t_h\delta_h)=r_g\afg(t_h1_{g\m})\omgh\delta_{gh}.\end{equation*}

By \cite[Theorem 2.4]{DES}  $R\star_{\af,\om} G$ is an associative ring with identity $1\delta_1.$ Moreover,  the map
$R\ni r\mapsto r\delta_1\in R\star_{\af,\om} G $ is a monomorphism of rings.  Thus, we consider $R$ as a subring of $R\star_{\af,\om} G,$ and the latter  as an $R$-$R$-bimodule via the actions
$$r(a_g\de_g)=ra_g\de_g,\,\,\,\,\text{and}\,\,\,\,\,(a_g\de_g)r=a_g\af_g(r1_{g\m})\de_g, \,\,\,r\in R,\, g\in G, \,a_g\in D_g.$$

\vu

The family of partial isomorphisms $\{\afg\colon D_{g\m}\to D_g\}_{g\in G}$ forms a  {\it partial action}
 (as defined
in \cite{DE}), which we also denote by $\af.$  Then  $\om=\{\omgh\}_{(g,h)\in G\times G}$ is called a {\it twisting} of $\af,$ and the
above twisted partial action will be simply denoted by $(\af, \om).$ Moreover, if $\om(g,h)=1_g1_{gh}$, for all $g,h\in G$,
we say that $\om$ is \emph{trivial}, and, in this case, the  partial crossed product  $R\star_{\af,\om} G$
coincides with the partial skew group ring $R\star_\af G$ as introduced in \cite{DE}.
The subring of {\it invariants} of $R$ under $\af$ was defined in \cite{DFP} as
\begin{equation*} \label{inva}R^\af=\{r\in R\,|\, \afg(r1_{g\m})=r1_g \,\,\forall g\in G\}.\end{equation*}

\vd

When $G$ is finite, the  extension $R\supseteq R^\af$ is  a {\it $\af$-partial Galois extension}    (see \cite{DFP}) if for some $m\in \N$ there exists a
subset $\{x_i, y_i\,|\, 1\le i\le m\}$ of $R,$ called a {\it partial Galois coordinate system} of the extension $R\supseteq R^\af,$ such that $\displaystyle\sum_{i=1}^{m}x_i\af_g(y_i 1_{g\m})=\delta_{1,g},\,\,\,\text{for all}\,\,\, g\in G.$

\begin{remark}\label{jjota}
It is shown in \cite[Theorem 4.1]{DFP} that  $R\supseteq R^\af$ is a partial Galois extension if and only if $R$ is a p.f.g $R^\af$-module and  the map
\begin{equation}\label{jota}
j\colon R\star_\af G \to {\rm End}_{R^\af}(R),\,\,\,\,\, j\left(\sum_{g\in G} r_g\de_g\right)(r)=\sum_{g\in G} r_g\afg(r1_{g\m}),
\end{equation} for each $r\in R$, is an isomorphism of $R$-modules and $R^{\af}$-algebras.
\end{remark}

\begin{defi}\cite[Definition 1.4]{DK} \label{defn-cochain}
Let  $T$ be a commutative ring or a monoid,  $n\in\N$ and $\af=(T_g, \af_g)_{g\in G}$ a unital partial action of $G$ on T.  An $n$-cochain of
$G$ with  values in $T$ is a map $f:G^n\to T,$ such that $f(g_1,\dots,g_n)\in \U(T1_{g_1}1_{g_1g_2}\dots1_{g_1g_2\dots g_n}),$
for all $n\in \N$. A $0$-cochain is an invertible element of $T$.
\end{defi}

The set of $n$-cochains $C^n(G,\af, T)$ is an abelian group under the point-wise multiplication. Its identity is
$(g_1, g_2\dots g_n)\mt 1_{g_1}1_{g_1g_2}\dots1_{g_1g_2\dots g_n}$ and the inverse of $f\in C^n(G,\af,T)$ is $f^{-1}(g_1,\dots,g_n)=f(g_1,\dots,g_n)^{-1}$,
where $f(g_1,\dots,g_n)^{-1}f(g_1,\dots,g_n)= 1_{g_1}1_{g_1g_2}\dots1_{g_1g_2\dots g_n}$.

\begin{prop}\cite[Proposition 1.5]{DK}\label{pcobh} Let $\delta^n\colon C^n(G,\af, T)\to C^{n+1}(G,\af, T) $ be the map defined by \begin{align*}\label{pcob}
(\delta^nf)(g_1,\dots,g_{n+1})=&\,\af_{g_1}\left(f(g_2,\dots,g_{n+1})1_{g_1^{-1}}\right)
\prod_{i=1}^nf(g_1,\dots , g_ig_{i+1}, \dots,g_{n+1})^{(-1)^i}\notag\\
&f(g_1,\dots,g_n)^{(-1)^{n+1}},
\end{align*} for any $f\in C^n(G,\af, T)$ and $g_1,\dots,g_{n+1}\in G,$ where  the inverse elements are taken in the corresponding ideals.
If $n=0$ and $t\in\U(T) $  we set $(\delta^0t)(g)=\af_g(1_{g^{-1}}t)t^{-1},$ for all $g\in G$. Then,
		$\delta^n$  is a homomorphism such that
	$
		(\delta^{n+1}\delta^nf)(g_1,g_2,\dots, g_{n+2})=1_{g_1}1_{g_1g_2}\dots1_{g_1g_2\dots g_{n+2}},
	$
for any $f\in C^n(G,\af,T)$.
	\end{prop}
\begin{defi}\label{defn-cohomology}
		The map $\delta^n$ is called  a coboundary homomorphism. We define the groups $Z^n(G,\af,T)=\ker{\delta^n}$, $B^n(G,\af,T)={\rm im}\,{\delta^{n-1}}$
and $H^n(G,\af,T)=\df{\ker{\delta^n}}{{\rm im}\,{\delta^{n-1}}}$ of  partial $n$-cocycles, $n$-co\-boun\-da\-ries and $n$-cohomologies of
$G$ with values in $T$, $n\ge 1,$ respectively. For $n=0$ we  set $H^0(G,\af,T)=Z^0(G,\af,T)=\ker{\delta^0}$. We say that $f$ and $f'$ in
$C^{n+1}(G,\af,T)$ are cohomologous if $f=(\de^{n}\rho) f',$ for some $\rho \in C^{n}(G,\af,T ).$
	\end{defi}

{\it From now on G will denote a finite group, $R$ a commutative ring  and $R\supseteq R^\af$ a partial Galois extension.}

\end{section}

\begin{section}{
The exact sequence\\ $0\to H^1(G,\af, R)\stackrel{\f_1}{\to} {\bf Pic}(R^\af)\stackrel{\f_2}{\to }{\bf PicS}(R)^{\af^*}\cap {\bf Pic}(R) \stackrel {\f_3}{\to }H^2(G,\af, R)$}\label{f_1f_2f_3}


\vu

  For any   left $R\star_\af G$-module $M$   we set
 $$ M^G=\{m\in M\,|\, (1_g\de_g)m=1_gm, \,\,\text{for all}\,\, g\in G\}.$$

\begin{lema}\label{isop} For every left $R\star_\af G$-module M the map $\mu\colon R\ot_{R^\af} M^G\to M$ given by $\mu(x\ot_{R^\af} m)=xm$
is an isomorphism of left $R\star_\af G$-modules, where  $R$ is a  left $R\star_\af G$-module via
$\label{triangle} (r_g\de_g)\triangleright r=r_g\afg(r1_{g\m}).$
\end{lema}
\p By (ii) of \cite[Theorem 4.1]{DFP} the map $\mu$ is an isomorphism of $R$-modules.    So, we only must check that $\mu$ is $R\star_\af G$-linear.
 For $g\in G, r\in R, r_g\in D_g$ and $m\in M^G$ we have
\begin{align*}
\mu((r_g\de_g)\triangleright (r\ot_{R^\af} m))&=r_g\afg(r1_{g\m})m
\\&=[r_g\afg(r1_{g\m}) (1_g\de_g)] m
\\& =[(r_g\de_g)(r\de_1)]m\\
&=(r_g\de_g)[(r\de_1)m]
\\
&=(r_g\de_g)(rm)
\\&=(r_g\de_g)\mu (r\ot_{R^\af} m),
\end{align*}
as desired.\cua
\vu

We denote by  ${\bf Pic}(R^\af)$ the Picard group of $R^\af.$  For  details on the Picard group of a commutative ring the  reader may consult  \cite{TF}  or \cite{L}.

We recall the group homomorphism $\varphi_1\colon H^1(G,\af, R) \to {\bf Pic}(R^\af)$ given in \cite{DPP1}.  First of all we have that
\begin{align*}
		Z^1(G,\af, R) &=\{f\in C^1 (G,\af,R)\mid f(gh)1_g=f(g) \, \af _g ( f(h) 1_{g^{-1}} ), \forall g,h\in G \}, \text{and} \\
	B^1(G,\af, R) &=\{f\in C^1 (G,\af,R)\mid f(g)= \af _g ( a1_{g^{-1}} )a\m, \,\,\text{for some}\,\, a\in \U(R) \}.\end{align*}

 Let $f\in Z^1(G,\af, R)$ and  let $\te_f\in {\rm End}_{R^\af}(R\star_\af G)$  be given by $\te_f(r_g\de_g)=r_gf(g)\de_g$ for all $r_g\in D_g, g\in G.$ Then $\te_f$ is an $R^\af$-algebra homomorphism, and we define a $R\star_\af G$-module $R_f$ by $R_f=R$ as sets and with the action of  $R\star_\af G$ induced by
\begin{equation}\label{rf}(r_g\de_g)\cdot  r= {\te_f(r_g\de_g) \triangleright r, } \,\,\,\text{for any}\,\, r\in R,\, g\in G. \end{equation} We know from \cite{DPP1} that
 $[R_f^G]\in {\bf Pic}(R^\af)$ and the map
\begin{equation*}\label{1homo}\varphi_1\colon H^1(G,\af, R) \ni {\rm cls}(f) \to [R_f^G]\in  {\bf Pic}(R^\af)\end{equation*}
is a well defined group homomorphism.

\vu

\begin{teo}\label{homp} The sequence  $0\to H^1(G,\af, R)\stackrel{ \varphi_1} \to {\bf Pic}(R^\af)$ is exact.
\end{teo}
\p   If $ \varphi_1({\rm cls}(f))=[R_f^G]=[R^\af],$ then by Lemma \ref{isop} there are $R\star_\af G$-module isomorphisms $R_f\cong R \ot_{R^\af} R_f^G \cong R.$
{ By \cite[Remark 2.8]{DPP1}  $f$ is normalized, i.e. $f(1)=1$.}  Let $w\in{\rm Hom}_{R\star_\af G}(R_f, R)$ be an isomorphism. Since $f(1)=1$ we have $w\in{\rm Hom}_{R^\af}(R, R)\cong R\star_\af G,$ by Remark \ref{jjota}.
 Hence,   there is an element $W=\sum\limits_{g\in G}r_g\de_g\in \U(R\star_\af G)$ such that  $w(x)=Wx,$ for all $x\in R.$ Recall that $R\star_\af G$ acts on $R_f$ via $V \cdot r=\te_f(V)r, \, r\in R, V\in R\star_\af G .$  Then, $$(VW)x=V(Wx)=V(w(x))=w(V\cdot x)=w(\te_f(V) x)=(W\te_f(V))x,$$ for all $x\in R_f=R,$ where the third equality above holds because   $w\in{\rm Hom}_{R\star_\af G}(R_f, R).$ Thus,
$W\te_f(V)=VW.$
 In particular, if $V=r\de_1=r,$ $\te_f(V)=r,$ and this gives $\sum\limits_{g\in G}r_g\afg(r1_{g\m})\de_g=\sum\limits_{g\in G}rr_g\de_g.$ Consequently,    $j(r_g\de_g)(r)=j(r_g\de_1)(r),$ for all $g\in G,$ where $j$ is the isomorphism defined in (\ref{jota}).
 Therefore, $r_g=0$ if $g\ne 1$ and $W=r_1\de_1\in \U(R).$ Now, taking $V=1_h\de_h$ we get $\te_f(V)=f(h)\de_h$ and
$$f(h)\de_h=(r\m_1\de_1)(1_h\de_h)(r_1\de_1)=\afh(r_11_{h\m})r\m_1\de_h,$$ which implies $f\in B^1(G,\af, R).$ \cua
\medskip

Let $\kappa$ be a commutative ring and $\Lambda$ be    a commutative  $\kappa$-algebra with identity. We recall from \cite{DPP1}  the following.
\begin{defi} A   non-necessarily central $\La$-$\La$-bimodule $P$ is called  $\kappa$-partially invertible if $P$ is a  central $\kappa$-bimodule and   as a    $\La$-$\La$-bimodule $P$ satisfies  the following two properties.
\begin{itemize}
\item P is  left and right  f.g.p. $\Lambda $-module,
\item  $\Lambda ^{op} \ni \lambda\mapsto  r_\lambda \in {\rm End}(_\Lambda P),
 r_\lambda(p)=p\lambda\,\, \text{  and }\,\,    \Lambda \ni \lambda\mapsto   l_\lambda \in {\rm End}(P_\Lambda), l_\lambda(p)=\lambda p,$ are $k$-algebra epimorphisms.
\end{itemize}

Let $[P]=\{ M\,|\, M\, \text{is a}\,\, \La\text{-}\La\text{-bimodule and} \,\, M\cong P\, \,\text{as} \,\,\La\text{-}\La\text{-bimodules}  \}.$ We denote by  ${\bf PicS}_\kappa(\La)$ the  set of classes $[P]$ of partially invertible $\La$-$\La$-bimodules with composition $[P][Q]=[P\ot_\La Q].$ If $\kappa=\Lambda,$ we set
${\bf PicS} (\La)  = {\bf PicS}_{\La}(\La).$
\end{defi}
  The following result is \cite[Proposition 3.5]{DPP1}. Since the proof presented there has several misprints, we give it below for reader's convenience. 

 \begin{prop}\label{ccom} The product $[P][Q]=[P\ot_\La Q]$ endows   ${\bf PicS}_{\kappa}(\La)$ with the structure of a  semigroup.
\end{prop}

 \p  We  shall show that  $[P\ot_\La Q]\in   {\bf PicS}_{k} (\La),$ for any
$[P], [Q]\in   {\bf PicS}_{\kappa }  (\La).$  Notice that $P\ot_\La Q$ is a
left f.g.p. $\La$-module. Indeed, there are free f.g. left  $\La$-modules $F_1,F_2$ and  left $\La$-modules $M_1,M_2$ such that
$P\oplus M_1=F_1,\, Q\oplus M_2=F_2.$ Now consider $M_1$ and $F_1$ as central $\La$-$\La$-bimodules, then by tensoring the two previous
equalities we see that there exists a left $\La$-module $M$ such that $(P\ot_\La Q)\oplus M  \cong  F_1\ot_\La F_2$   as left $\La$-modules, and the assertion follows. In a similar way, one can show that  $P\ot_\La Q$ is a
right f.g.p. $\La$-module.

\vu

By assumption there are   $\kappa $-algebra epimorphisms $r^1 \colon \Lambda \to  {\rm End}(_\Lambda P)$ and  $r^2\colon \Lambda \to  {\rm End}(_\Lambda Q),$ given by right multiplications.   It follows, using the third item of Lemma~\ref{ht},  that $$r^1\ot r^2\colon \Lambda\ot_\Lambda
\Lambda\to {\rm End} (_{\La} P) \ot_{\La } {\rm End} ( _{\La}  Q) \cong {\rm End}_{\La}(P\ot_\La Q)$$ is a   $\kappa $-algebra epimorphism. Since $\Lambda\ni  \lambda  \mt 1_\Lambda\ot_\Lambda\lb\in
\Lambda\ot_\Lambda \Lambda$ is a   $\kappa $-algebra isomorphism, we conclude that  $ \Lambda \ni \lambda\mapsto
r^1_{1_\Lambda }\ot r^2_\lambda\in {\rm End}_{\La}(P\ot_\La Q)$ is a  $\kappa $-algebra epimorphism.   Now for any $p\ot_\La q\in P\ot_\La Q,$ we have that $r^1_{1_\Lambda }\ot r^2_\lambda(p\ot_\La q)=p\ot_\La q\lambda,$   and  the assertion follows. In an analogous way we  obtain that the left multiplication is a  $\kappa $-algebra epimorphism.\cua 

\begin{defi} \cite[Definition 3.1]{DPP1} We say that a f.g.p. central $R$-$R$-bimodule  P has rank less than or equal to one, if for any  $\mfp\in {\rm Spec}(R)$ one has $P_\mfp=0$ or $P_\mfp\cong R_\mfp$ as $R_\mfp$-modules. In this case we write ${\rm rk}_{R}(P)\le 1.$
\end{defi}
We remind a characterization of ${\bf PicS} (R).$
\begin{prop} \label{picinv}  For any commutative ring R we have
\begin{itemize}
\item  \cite[Definition 3.1, Remark 3.4]{DPP1}
$${\bf PicS} (R)=\{[E]\,|\, E\, \text{is a f.g.p. $R$-module and}\,\, {\rm rk}(E)\le 1 \}.$$

\vu

\item   \cite[Proposition 3.8, arXiv version]{DPP1}.

\vu

The set  ${\bf PicS} (R),$ with binary operation induced by the tensor product, is a commutative inverse  monoid with 0. Moreover,  the inverse $[E]^*$ of $[E]$ is $[E^*],$ where $E^*={\rm hom}_R(E,R)$ is the dual of $E.$

\vu

\item \cite[Formula (3.1)]{DPP1}
\begin{equation*}  {\bf PicS} (R)\cong\bigcup\limits_{e\in {\bf I_p}(R)}{\bf Pic} (Re),
\end{equation*} where ${\bf I_p}(R)$ denotes the semilattice of the idempotents of $R$ with respect to the product.

\end{itemize}
\end{prop}

\vd

{ We recall also the next.} 
\begin{lema} \cite[Lemma 3.12, arXiv version]{DPP1}\label{iguald}  Let $g\in G$ and $[M]\in {\bf PicS}(R).$  If $1_gm=m,$ for all $m\in M$ and $M_\mfp\cong (D_g)_\mfp$ as $R_\mfp$-modules, for all $\mfp\in{\rm Spec}(R),$ then   $[M]\in {\bf Pic}(D_g).$
\end{lema}

\vd

For the sequel, we proceed by recalling  a partial action of $G$ on  PicS($R$) given in \cite{DPP1}. Let $\af=(D_g, \afg)_{g\in G}$ be a partial action of $G$ on $R.$
 Then, for any $y\in R$  we have
\begin{equation}\label{prodp}\afg(\afh(y1_{h\m})1_{g\m})=\af_{gh}(y1_{(gh)\m})1_g, \,\,\,\text{ for all}\,\, g,h\in G. \end{equation}

For any $D_{g\m}$-module $E,$ we denote by $E_g$  the $R$-module $E$ via the  action
\begin{equation}\label{bulletaction}
r\bullet x=\af_{g\m}(r1_g)x,
\end{equation}
for any $r\in R,\, x\in E.$
\begin{lema}\label{gaction} We have the following.
\begin{itemize}
\item \cite[Lemma 3.6, vi)]{DPP1}   For any $[M]\in {\bf Pic}(D_{g\m}),$ $[M_g]\in{\bf Pic}(D_{g}). $

\vu

\item \cite[ Lemma 3.7]{DPP1} $X_g=[D_g]{\bf PicS}(R)=\{[E]\in {\bf PicS} (R)\,|\, E=1_gE\}.$

\vu

\item \cite[ Theorem 3.8]{DPP1} The family $\af^*=( X_g, \afg^*)_{g\in G},$ where $$\afg^*\colon X_{g\m}\ni [E]\mapsto [E_g]\in X_{g},$$ gives a partial action of G on the inverse semigroup ${\bf PicS} (R).$   Moreover, $\U(X_g)={\bf Pic}(D_g).$ 

\vu

\item \cite[Proposition 4.5, arXiv version]{DPP1}  $${\bf PicS}(R)^{\af^*}=\{[E]\in {\bf PicS}(R)\,|\, (E\ot D_{g\m})_g\cong E\ot D_g, \,\,\text{for all} \,g\in G\} $$  has a 0 and is a commutative  inverse submonoid of ${\bf PicS}(R).$ Moreover, \begin{equation}\label{esg}[(E^*)_g]=[E_g]^*=[(E_g)^*],\end{equation}
for any $g\in G.$
\end{itemize}
\end{lema}

\vu

\begin{prop}\cite[ Proposition 4.2]{DPP1} The map $$\varphi_2\colon  {\bf Pic}(R^\af)\ni [E]\mapsto[ R\ot_{R ^\af} E]\in {\bf PicS}(R)^{\af^*}\cap {\bf Pic}(R),$$

is a group homomorphism.
\end{prop}

\vu

\begin{teo} The sequence $H^1(G,\af, R)\stackrel{\varphi_1}{\to} {\bf Pic}(R^\af)\stackrel{\varphi_2}{\to }{\bf PicS}(R)^{\af^*}\cap {\bf Pic}(R)$ is exact.
\end{teo}
	\p Let $f\in Z^1(G,\af, R).$ Then $\varphi_2\varphi_1({\rm cls}(f))= [R\ot_{R^\af} R^G_f]$ and
$  R\ot_{R^\af} R^G_f \cong R_f,$ the latter being  the isomorphism of $R\star_\af G$-modules given by Lemma \ref{isop}. Thus,  $ R\ot_{R^\af} R^G_f \cong R$ as $R$-modules, and ${\rm im}\,\varphi_1\subseteq \ker \varphi_2.$

\vu

Let $[M]\in {\bf Pic}(R^\af) $ and assume $[M]\in \ker \varphi_2.$ Then, there is an $R$-module isomorphism $\xi\colon R\ot_{R^\af} M\to R.$ Observe that $\xi(D_g\ot _{R^\af} M)=1_g\xi(R\ot _{R^\af} M)=1_gR=D_g.$ Since $R$ is an $R\star_\af G$-module, with action $\triangleright$ defined in  Lemma~\ref{triangle},  $R\ot _{R^\af} M$ is also an $R\star_\af G$-module by  letting $R\star_\af G$ act on $R.$ We identify $R\star_\af G$ with ${\rm End }_{R^\af}(R)$ via (\ref{jota}) and define a map $\theta\in {\rm End }_{R^\af}(R\star_\af G)$ by
\begin{equation}\label{dete}\te(W)( x)=\xi(W\triangleright\xi\m(x)),\,\,\,\text{ for all}\,\,\, W\in R\star_\af G\,\,\, \text{and} \,\,\, x\in R. \end{equation}
Notice that $\te$ is an  $R^\af$-algebra homomorphism and $\te(r)( x)=rx$ for all $r,x\in R.$

Let $f\colon G\ni g\mapsto \te(1_g\de_g) (1)\in R,$ where $\te(1_g\de_g)$ is viewed as an element of ${\rm End }_{R^\af}(R).$
We shall prove that $f\in Z^1(G,\af, R).$

\vu

First of all, since $1_g\de_g  \triangleright\xi\m(1)\in D_g\ot _{R^\af} M$  we see that $f(g)\in D_g,$ for all $g\in G.$  On the other hand,
\begin{align*}
\te(1_g\de_g)( r)=\xi(1_g\de_g \triangleright\xi\m(r))=\xi(\afg(r1_{g\m})\de_g \triangleright\xi\m(1))=\afg(r1_{g\m})\xi(1_g\de_g \triangleright\xi\m(1)),
\end{align*}
and so
\begin{equation}\label{iguall}\te(1_g\de_g)( r)=\afg(r1_{g\m})f(g), \,\,\,\,\, \text{for each}\,\,\, g\in G\,\,\, \text{and}\,\,\, r\in R.\end{equation}
 In particular, for  $r_g=\te(1_{g\m}\de_{g\m})( 1_g),$ we have
\begin{align*}
\afg(r_g1_{g\m})f(g)&=\te(1_g\de_g) (\te(1_{g\m}\de_{g\m})(1_g))\\
&=[\te(1_g\de_g)\te(1_{g\m}\de_{g\m}) ]( 1_g)\\
&=\te(1_g\de_1)( 1_g)
\\&=1_g,
\end{align*}
ensuring that $f(g)\in \U(D_g),$ for all $g\in G.$ Thus,  $f\in C^1(G,\af,R).$

\vu

Now we show that $f$ satisfies the $(1,\af)$-cocycle identity.
By  (\ref{iguall}) we have $\te(1_g\de_g)( r)=( f(g)\de_g) (r),$ for any $r\in R.$ Then,  $\te(1_g\de_g)=f(g)\de_g$ as maps defined on $R,$ and
\begin{align*}
f(g)\afg(f(h)1_{g\m})\de_{gh}=(f(g)\de_g)(f(h)\de_h)=\te(1_g\de_g)\te(1_h\de_h)=1_g\te(1_{gh}\de_{gh})=1_gf(gh)\de_{gh},
\end{align*}
showing that $f\in Z^1(G,\af, R).$
To end the proof we must show that there is an $R^\af$-module isomorphism from $M$ to $R^G_f.$ Consider the map $\xi\m\colon R_f\to R\ot  _{R^\af} M,$ then for each $z \in R_f$ we have
\begin{align*}
\xi\m((r_g\de_g)\cdot z)\stackrel{(\ref{rf})}{=}&\,\xi\m(r_g\afg(z1_{g\m}) { f(g) } )\stackrel{(\ref{iguall})}{=}  \xi\m(\te(r_g\de_g)(z))\stackrel{(\ref{dete})}{=}(r_g\de_g)\triangleright \xi\m(z).
\end{align*}
Hence, there are isomorphisms of  $R\star_\af G$-modules $R\ot  _{R^\af} R^G_f\cong R_f\cong R\ot  _{R^\af} M,$ which leads to $R^G_f\cong M$ as $R^\af$-modules, thanks to the second item of Lemma \ref{proje}. Therefore, $[M]=\f_1({\rm cls}(f))$.\cua

\vt

We recall from \cite{DPP1} the   homomorphism $\varphi_3\colon {\bf PicS}(R)^{\af^*}\cap {\bf Pic}(R)\to H^2(G,\af,R).$   First of all we have that $B^2 (G,\af,R)$ is the group \begin{align*}
&\left\{ w \in  C^2 (G,\af, R) \mid \exists  f \in  C^1 (G,\af, R),\text{with}\,\,w(g,h) = \af _g ( f(h) 1_{g^{-1}} )f(g) f(gh)\m
\right\}.\end{align*}
	 and $Z^2(G,\af,R)$ is
	\begin{align*}
		&\{w\in C^2 (G,\af, R) \mid \af _g ( w(h,l) 1_{g^{-1}} ) \; w(g,hl) = w(gh,l)  \, w(g,h),\,\,\forall g,h,l \in G \}.\end{align*}

Now, let $[E]$ be an element of ${\bf PicS}(R)^{\af^*}\cap {\bf Pic}(R).$ Then, there is a family of $R$-module isomorphisms
\begin{equation}\label{psiii}\{\psi_g\colon ED_g\to (ED_{g\m})_g\}_{g\in G},\,\, \text{with}\,\, \psi_g(rx)=\af_{g\m}(r1_g)\psi_g(x), \,\end{equation}
{ for all $r\in R, x\in ED_g, g\in G, $\footnote{ Since $E\ot D_g\ni x\ot  d\mt xd \in ED_g$ is a $R$-module isomorphism, we identify $E\ot D_g$ with $ED_g,$ for any $g\in G.$  }
and the maps, $\psi\m_g\colon (ED_{g\m})_g\to ED_g, g\in G$, satisfy
\begin{equation}\label{psii}\psi\m_g(rx)=\afg(r1_{g\m})\psi\m_g(x),\,\,\text{ for all}\,\, g\in G\,\,,r\in R, \,\,x=x_{g\m}\in ED_{g\m}.
\end{equation} }

\vu

Hence, $\psi_{(gh)\m}\psi\m_{h\m}\psi\m_{g\m}\colon ED_gD_{gh}\to  ED_gD_{gh}$ belongs to $\U({\rm End}_{D_gD_{gh}}(ED_gD_{gh})).$
Since $[E]\in  {\bf Pic}(R)$,\, it follows that $[ED_gD_{gh}]\in {\bf Pic}(D_gD_{gh}),$ and consequently we have ${\rm End}_{D_gD_{gh}}(ED_gD_{gh})\cong D_gD_{gh}.$ Therefore,   there exists $\om_{g,h}\in \U(D_gD_{gh})$ such that $$\psi_{(gh)\m}\psi\m_{h\m}\psi\m_{g\m}(x)=\om_{g,h}x,\,\,\text{ for all}\,\, x\in ED_gD_{gh}.$$
Summarizing, for an element $[E]\in {\bf PicS}(R)^{\af^*}\cap {\bf Pic}(R) $ we have found a map $\om_{[E]}=\om\colon G\times G\ni (g,h)\mapsto \om_{g,h}\in \U(D_gD_{gh})\subseteq R.$ Moreover, $\om \in Z^2(G,\af,R)$ and   by  \cite[Claim 4.3]{DPP1}  ${\rm cls(\om)}$ does not depend on the choice of the family  of isomorphisms given by \eqref{psiii}.

\begin{prop}\cite[Claim 4.4]{DPP1}  $\varphi_3\colon\! {\bf PicS}(R)^{\af^*}\cap {\bf Pic}(R)\ni [E] \mapsto {\rm cls(\om)}\!\in\! H^2(G,\af, R)$ is a  group homomorphism.\end{prop}

We need  the following.
\begin{lema}\label{cond} Let $[E]\in  {\bf PicS}(R)^{\af^*}$  and  $\{\psi_g\colon ED_g\to (ED_{g\m})_g\}_{g\in G}$ be a family of R-module isomorphisms with $\psi_1={\rm id}_E.$ Suppose that   $c_g$ are fixed elements  in $D_g,\, g\in G$, with $c_1 =1.$ Then the action $$\left(\sum\limits_{g\in G}a_g\de_g\right) x=\sum\limits_{g\in G}a_gc_g\psi_{g\m}(x1_{g\m})$$  endows the $R$-module { $E$ } with an $R\star_{\af,\om}G$-module structure, provided that
\begin{equation}\label{condit}\om_{g,h}c_{gh}\psi_{(gh)\m}(x)=c_g\afg(c_h1_{g\m})\psi_{g\m}\psi_{h\m}(x)\end{equation} holds for all $g,h\in G$ and  $x\in ED_{h\m}D_{(gh)\m}.$
\end{lema}
\p Since $\psi_1={\rm id}_E$ and the action is linear we only need to check  the equality $[(a_g\de_g)(a_h\de_h)]x=(a_g\de_g)[(a_h\de_h) x]$ for all $g,h\in G$ and $x\in E.$ From (\ref{psii}) we obtain that
\begin{equation}\label{palfa}
\psi_{g\m}(\psi_{h\m}(x1_{h\m})1_{g\m})=\psi_{g\m}(\psi_{h\m}(x1_{(gh)\m}1_{h\m}))\,\,\,\,\text{for all}\,\, g,h\in G, x\in E.
\end{equation}
Then
\begin{align*}[(a_g\de_g)(a_h\de_h)] x&=a_g\afg (a_h1_{g\m})\om_{g,h}c_{gh}\psi_{(gh)\m}(x1_{(gh)\m})\\
&=a_g\afg (a_h1_{g\m})\om_{g,h}1_g 1_{gh}c_{gh}\psi_{(gh)\m}(x1_{(gh)\m})\\
&=a_g\afg (a_h1_{g\m})\om_{g,h}c_{gh}\psi_{(gh)\m}(x1_{(gh)\m})\af_{gh}(1_{h\m}1_{(gh)\m})\\
&=a_g\afg (a_h1_{g\m})\om_{g,h}c_{gh}\psi_{(gh)\m}(x1_{(gh)\m}1_{h\m}).
\end{align*}
On the other hand
\begin{align*}(a_g\de_g) [(a_h\de_h) x]&=(a_g\de_g)(a_hc_h\psi_{h\m}(x1_{h\m}))\\
&=a_gc_g\psi_{g\m}(a_hc_h\psi_{h\m}(x1_{h\m})1_{g\m})\\
&=a_gc_g\afg (a_hc_h1_{g\m})\psi_{g\m}(\psi_{h\m}(x1_{h\m})1_{g\m})\\
&\stackrel{(\ref{palfa})}=a_g\afg (a_h1_{g\m})c_g\afg(c_{h}1_{g\m})\psi_{g\m}\psi_{h\m}(x1_{(gh)\m}1_{h\m}).
\end{align*}
Using (\ref{condit}) the  Lemma follows.\cua

\vu

\begin{teo} The sequence ${\bf Pic}(R^\af)\stackrel{\varphi_2}{\to }{\bf PicS}(R)^{\af^*}\cap {\bf Pic}(R)\stackrel {\varphi_3}{\to }H^2(G,\af, R)$ is exact.
\end{teo}
\p  If    $[E]\in {\bf Pic}(R^\af) ,$ then   $\varphi_2([E])=[R\ot_{R^\af} E]$ by the definition of $\f _2.$
 There  are isomorphisms $(R\ot_{R^\af} E)\ot D_g\cong  D_g\ot_{R^\af} E$ and $((R\ot_{R^\af} E)\ot D_{g\m})_g\cong ( D_{g\m}\ot_{R^\af}  E)_g,$ of $R$-modules.
 Notice that $D_g\ot_{R^\af}  E\cong ( D_{g\m}\ot_{R^\af} E)_g$ as $R$-modules via the map
$d\ot_{R^\af}  x\stackrel{\psi_g}{\mapsto} \af_{g\m}(d)\ot_{R^\af}, $ where $ x, \,d\in { D_{g}}, x\in E, g\in G.$ Then $\psi_{(gh)\m}\psi\m_{h\m}\psi\m_{g\m}(d x)=1_g1_{gh}(d x),$ for all  $d\in D_{g}D_{gh} $ and any $x\in E, g,h\in G.$ Therefore ${\rm im}\,\varphi_2\subseteq \ker \varphi_3.$

\vu

Let $[V]\in \ker \varphi_3.$ Then $ \varphi_3([V])={\rm cls(\om)},$ with $\om \in B^2(G,\af,R).$  There are $R$-module isomorphisms $\psi_g\colon V D_g\to  (VD_{g\m})_g$ and a map $u\colon G\ni g\to u_g\in \U(D_g)\subseteq R$ such that $$\psi_{(gh)\m}\psi\m_{h\m}\psi\m_{g\m}(x)=\om_{g,h}x=u_g\afg(u_h1_{g\m})u\m_{gh}x, \,\,\,\,x\in { V}D_gD_{gh},\,\,\, g,h\in G.$$ Since cls$(\om)$ does not depend on the choice of the isomorphisms, we can choose $\psi_1= {{\rm id}_V. } $  Lemma \ref{cond} and (\ref{psiii}) imply that
we can view $V$ as an $R\star_\af G$-module via the action
$$(a_g\de_g) x=a_gu_g\psi_{g\m}(x1_{g\m}),$$
for all $a_g\in D_g, x\in V$ and $g\in G$. By Lemma \ref{isop}  there is an $R\star_\af G$-module isomorphism $V\cong R\ot_{R^\af} V^G.$ Notice that $V$ is a  f.g.p. $R^\af$-module, therefore  by the first item of Lemma \ref{proje}, $V^G$ is a   f.g.p. $R^\af$-module. Hence, it only remains to prove that $[V^G]\in {\bf Pic}(R^\af), $ in order to show that $V\in\text{im}\,\varphi_2$.
Let $\mfp\in {\rm Spec}(R^\af).$ The fact that  $R$ is a f.g.p.  $R^\af$-module and  \cite[Proposition 20.6]{LR} imply that  $R_\mfp$ is  a semi-local ring, and, consequently ${\bf Pic}(R_\mfp)=0$ (see \cite[Example 2.22 (D)]{L}) which leads to $V_\mfp \cong R_\mfp$ as $R_\mfp$-modules. Thus, $R_\mfp\ot_{(R^\af)_\mfp}(R^\af)_\mfp\cong R_\mfp\cong V_\mfp \cong R_\mfp\ot_{(R^\af)_\mfp}(V^G)_\mfp$ as $(R^\af)_\mfp$-modules. The required follows now from the the second item of Lemma 2.1.\cua

\end{section}

\section{The homomorphism $\f_4$ and  exactness at $H^2(G,\af, R)$}


\vu

\begin{prop} \cite[Theorem 5.7]{DPP1} The map $ \varphi_4\colon H^2(G,\af, R)\ni {\rm cls}(\om)\mapsto [R\star_{\af,\om} G]\in B(R/R^\af)$
 is a group homomorphism.\end{prop}

\begin{teo} \label{fi3fi4}The sequence ${\bf PicS}(R)^{\af^*}\cap {\bf Pic}(R)\stackrel{\varphi_3}{\to}H^2(G,\af, R)\stackrel{\varphi_4}{\to} B(R/R^\af)$ is exact.
\end{teo}
\p Let $[E]$ be an element of ${\bf PicS}(R)^{\af^*}\cap {\bf Pic}(R)$ and write $\varphi_3([E])={\rm cls(\om)}.$ Let $\{\psi_g\colon ED_g\to (ED_{g\m})_g\mid g\in G\},$ (where $\psi_1={\rm id}_E$) be a family of $R$-module isomorphisms determining $\om.$ That is, $\om_{g,h}x=\psi_{(gh)\m}\psi\m_{h\m}\psi\m_{g\m}x,$ for all $g,h\in G$ and all $x\in ED_gD_{gh}.$ Then, it is immediate to see that (\ref{condit}), with  $c_g =1_g$ and $\omega$  replaced by $\omega \m ,$ is satisfied, and by Lemma \ref{cond}  we may consider $E$  as an $R\star_{\af, \om\m}G$-module, where the  action of $R\star_{\af, \om\m}G$ on $E$ is induced by
$$(a_g\de_g)x=a_g\psi_{g\m}(x1_{g\m}), \, x\in E, g\in G, a_g\in D_g.$$
(Notice that with this new action the structure of $E$ as an $R$-module does not change.)
 Hence, there is an $R^\af$-algebra homomorphism $\zeta\colon R\star_{\af, \om\m}G\to {\rm End}_{R^\af}(E),\, \zeta(W)(x)=W x,\, $ for all $W\in R\star_{\af, \om\m}G$ and $ x\in E. $ We shall  prove that $\zeta$ is an isomorphism.
 Set $C=C_{{\rm End}_{R^\af}(E)}(R\star_{\af, \om\m}G).$ By \cite[Theorem 2.4.3]{DI} there is an $R^\af$-algebra isomorphism
\begin{equation}\label{multipl}R\star_{\af, \om\m}G \; \ot_{R^\af} C\cong {\rm End}_{R^\af}(E),\end{equation}
given by multiplication.  Since $R\star_{\af, \om\m}G$ and $ {\rm End}_{R^\af}(E)$ are f.g.p. $R^\af$-modules and $(R\star_{\af, \om\m}G)_\mfp\ne 0$ for all $\mfp \in {\rm Spec}(R^\af),$ the first item of Lemma \ref{proje} and (\ref{multipl}) imply that $C$ is a f.g.p. $R^\af$-module. Now, for any $\mfp \in {\rm Spec}(R^\af),$ the ring $R_\mfp$ is semilocal and $E_\mfp\cong R_\mfp$ as $ R_\mfp$-modules, because $[E]\in{\bf Pic}(R).$
 Thus, as $R^\af_{\mfp}$-modules  ${\rm End}_{R^\af}(E)_\mfp\cong {\rm End}_{R^\af}(R)_\mfp\cong (R\star_{\af, \om\m}G)_\mfp,$ where the latter isomorphism is $ j_\mfp$ induced by (\ref{jota}). Then, there exists $n_\mfp \in \N$  such that ${\rm End}_{R^\af}(E)_\mfp \cong (R\star_{\af, \om\m}G)_\mfp \cong (R^\af)_\mfp^{n_\mfp},$ and we conclude from  (\ref{multipl}) that ${\rm rk}_{R^\af}(C)=1.$ Then there exists $u_\mfp\in C_\mfp$ such that $C_\mfp=(R^\af)_\mfp u_\mfp$ and  localizing  (\ref{multipl}) we have
 $${\rm End}_{R^\af}(E)_\mfp=(R\star_{\af, \om\m}G)_\mfp C_\mfp=(R\star_{\af, \om\m}G)_\mfp(R^\af)_\mfp u_\mfp=(R\star_{\af, \om\m}G)_\mfp u_\mfp,$$
 from which we see that $u_\mfp$ is a unit  in ${\rm End}_{R^\af}(E)_\mfp$ and ${\rm End}_{R^\af}(E)_\mfp=(R\star_{\af, \om\m}G)_\mfp.$ Since $\mfp$ is arbitrary and ${\rm End}_{R^\af}(E)\supseteq R\star_{\af, \om\m}G,$ then  ${\rm End}_{R^\af}(E)=R\star_{\af, \om\m}G.$ Furthermore, $E$ is a faithfully projective $R^\af$-module and we get $[R\star_{\af, \om\m}G]=[{\rm End}_{R^\af}(E)]=[R^\af].$ Then ${\rm cls}(\om\m)\in \ker\varphi_4,$ and, consequently ${\rm cls}(\om)\in \ker\varphi_4.$  Hence ${\rm im}\,\varphi_3\subseteq \ker \varphi_4.$

To prove the converse inclusion take  ${\rm cls}(\om)\in \ker \varphi_4.$ Then ${\rm cls}(\om\m)\in \ker \varphi_4$ and thus $[R\star_{\af, \om\m}G]=[R^\af].$ By \cite[Prop. 5.3]{AG} there exists a  faithfully projective $R^\af$-module $V$ and a $R^\af$-algebra isomorphism  $\Upsilon\colon R\star_{\af, \om\m}G\to {\rm End}_{R^\af}(V),$ and we may view $V$ as a (faithful) $R\star_{\af, \om\m}G$-module (and, consequently, as a faithful  and central $R$-bimodule) via
$W\diamond v=\Upsilon(W)(v).$

\vu

Recall that $R$ is a maximal commutative $R^{\af }$-subalgebra in $R\star_{\af, \om\m}G$ (see \cite[Lemma 2.1 vi), Prop. 3.2]{PS}). We may consider a copy of $R$ in  ${\rm End}_{R^\af}(V)$ via the identification  $R=\Upsilon(R\delta_1).$ Then ${\rm End}_{R^\af}(V)$ is an Azumaya $R^\af$-algebra and $R=\Upsilon(R\delta_1)$ is a maximal commutative $R^\af$-subalgebra of ${\rm End}_{R^\af}(V).$ Since $R\supseteq R^\af$ is separable, \cite[Theorem 5.6]{AG} implies that ${\rm End}_{R^\af}(V)$ is f.g.p. $R$-module. Moreover, by \cite[Theorem 2]{AZ} $V$ is a right f.g.p. ${\rm End}_{R^\af}(V)$-module, and we see that $V$ is a f.g.p. $R$-module.  Furthermore, since ${\rm End}_{R}(V)=C_{{\rm End}_{R^\af}(V)}(R)=C_{R\star_{\af, \om\m}G}(R)=R,$  we conclude that $[V]\in {\bf Pic}(R).$

We show that  $[V]$ is an element of ${\bf PicS}(R)^{\af^*}.$ First notice  that the map
$$V\times D_g\ni (v,d)\mt (1_{g\m}\de_{g\m}\di v)\ot \af_{g\m}(d)\in (V\ot D_{g\m})_g, \,g\in G$$ is $R$-balanced. Indeed,
\begin{align*}
(v\di r, d)&\mt (1_{g\m}\de_{g\m})\di (v\di r)\ot \af_{g\m}(d)\\
&=  ((1_{g\m}\de_{g\m}) (r \delta  _1))\di  v\ot \af_{g\m}(d)\\
&=\af_{g\m}(r1_{g})\de_{g\m}\di v\ot \af_{g\m}(d)\\
&=((\af_{g\m}(r1_{g})\de_1) (1_{g\m}\de_{g\m}))\di v)\ot \af_{g\m}(d)\\
&=\af_{g\m}(r1_{g})\de_1 \di  (1_{g\m}\de_{g\m} \di v)\ot \af_{g\m}(d)\\
&=(1_{g\m}\de_{g\m}  \di v )\di(\af_{g\m}(r1_{g})\de_1) \ot \af_{g\m}(d)\\
&=(1_{g\m}\de_{g\m}\di v)\ot \af_{g\m}(rd),\
\end{align*} and hence we obtain a well defined map $\psi_g\colon V\otimes D_g\to  (V\ot  D_{g\m})_g$ such that
$$\psi_g(v\ot d)= (1_{g\m}\de_{g\m})\di v\ot \af_{g\m}(d),\,\,\,\,\text{for all}\,\, \,\,v\in V, d\in D_g. $$
Now given $r\in R,$ we have as above
 \begin{align*}
\psi_g(r\di (v\ot d))&=\psi_g((r\di v)\ot  d) \\
&= [(1_{g\m}\de_{g\m})\di (r\di v)]\ot  \af_{g\m}(d) \\& =  [ (\af_{g\m}(r1_{g}) \de_{g\m})   \di v ] \ot  \af_{g\m}(d)\\
&=[(\af_{g\m}(r1_{g})\de_1)\di  ( ( 1_{g\m}\de_{g\m})\di v ) ]\ot  \af_{g\m}(d) \\&=\af_{g\m}(r1_{g})\di\psi_g(v\ot d),
\end{align*}
and $\psi_g$ is $R$-linear. Analogously the map $\lb_g\colon (V\ot D_{g\m})_g\to V\otimes D_g,$ satisfying $\lb_g(v\ot d)=(\om_{g,g\m}\de_g)\di v\ot_R \afg(d),$ for all $v\in V,d\in D_{g\m}$ is well defined and using that $\alpha _{g\m }(\om_{g,g\m}) = \om _{g\m ,g }  $ we obtain
$$(\psi_g\circ \lb_g)(v\ot d)=((1_{g\m}\de_{g\m})(\om_{g,g\m}\de_g))\di v\ot d=(1_{g\m} \delta _1 \di v)\ot d= v\ot  1_{g\m}d= v\ot d,$$
where the last equality follows because $d\in D_{g\m}.$ Hence $\psi_g\circ \lb_g={\rm id}_{ (V\ot D_{g\m})_g},$ and analogously $\lb_g\circ \psi_g={\rm id}_{ V\ot D_{g}}.$ Therefore $\psi_g$ is a $R$-module isomorphism for which $\psi\m_g=\lb_g$ and, consequently,   $[V]\in {\bf PicS}(R)^{\af^*}.$   Finally, for $g,h\in G,\,v\in V$ and $d\in D_gD_{gh}$ we see that
 \begin{align*}\psi_{(gh)\m}\psi\m_{h\m}\psi\m_{g\m}(v\ot d)&=[(1_{gh}\de_{gh})(\om_{h\m,h}\de_{h\m})(\om_{g\m,g}\de_{g\m})]\di v\ot_R d\\
&=[(\af_{gh}(\om_{h\m,h}1_{(gh)\m})\om\m_{gh,h\m}\de_g)(\om_{g\m,g}\de_{g\m})]\di v\ot d\\
&\stackrel{(v)}{=}[(\om_{g,h}\de_g)(\om_{g\m,g}\de_{g\m})\di v]\ot d\\
&=[(\om_{g,h}\af_g(\om_{g\m,g})\om\m_{g,g\m}\de_{1})\di v]\ot d\\
&\stackrel{(\ref{afom})}{=}[(\om_{g,h}\de_1)\di v]\ot d.
\end{align*}
This yields that  $\psi_3([V])={\rm cls}(\om)$ and ${\rm im}\,\varphi_3\supseteq \ker \varphi_4.$\cua

\section{The homomorphisms $\f_5$ and exactness at the Brauer group}


We  proceed by recalling the construction of   the group homomorphism $\f_5\colon B(R/R^\af)\to H^1(G,\af^*,{\bf PicS}(R)\!)$  (see \cite[Section 5.1]{DPP1}  for details). Let $\af^*=(\af^*_g, X_g)_{g\in G}$ be the partial action of $G$ on ${\bf PicS}(R)$ given in the third item of Lemma  \ref{gaction}. Then
$B^1(G,\af^*, {\bf PicS}(R))$ is  the group $$\{ f\in C^1(G,\af ^*,{\bf PicS}(R)) \mid f(g) = \af^*_g([P] [D_{g^{-1}}])  [P^{*}], \, \mbox{for some}\,\, [P]\in  {\bf Pic}(R)  \}$$
and $Z^1(G,\af^*, {\bf PicS}(R))$ is given by
 \begin{align*}
	 \{f\in C^1 (G,\af ^*,{\bf PicS}(R))\mid f(gh)[D_g]=f(g) \, \af^* _g ( f(h) [D_{g^{-1}}] ),\,\, \forall g,h\in G \}.\end{align*}

 Let $f\in C^1 (G,\af ^*,{\bf PicS}(R)), \, g\in G$ and $\mfp \in {\rm Spec}(R).$ For simplicity of notation we shall write $f(g)_\mfp$ for the representative of the class $f(g)$ localized at $\mfp.$

\medskip

Let $[A]\in B(R/R^\af).$ By \cite[Theorem 5.7]{AG}  we may assume that $A$ contains $R$ as a maximal commutative subalgebra.
 Define the  left $R\ot_{R^\af}A^{op}$-module  $_g(1_{g\m}A) = 1_{g\m} A $   via
\begin{equation} \label{ngaction} (r\ot_{R^\af} a)\bullet a_g=\af_{g\m}(r1_g)a_ga.\end{equation}
It follows from  (\ref{ngaction}) that  { $_g(1_{g\m}A)$} is an object in the category ${}_{R\ot_{R^\af} A^{op}}\,{\rm Mod}$   and there is a unique central   $D_g$-module $M^{(g)}$ such that
\begin{equation}\label{unico}_g(1_{g\m} A)\cong M^{(g)} \ot  A,\,\,\text{as}\,\, (R\ot_{R^\af} A^{op})\text{-modules},\end{equation}
 where $M^{(g)}$ is considered as an  $R$-module via the map $r\mt r1_g,$  and  $M^{(g)} \ot  A$ is an $(R\ot_{R^\af} A^{op})$-module by means of  $$(r\ot_{R^\af} a) ( x \ot  b)  =  rx \ot  ba ,$$
$r\in R, a  \in  A^{op}  ,  x\in M^{(g)}, b\in A.$ 
Set
\begin{equation}\label{fa} f_A\colon G\ni g\mapsto [M^{(g)}]\in {\bf PicS}(R).\end{equation}

\begin{prop}\cite[Theorem 5.9]{DPP1}\label{fi5}
The map $$\varphi_5 \colon B(R/R^\af) \ni [A]\mt {\rm cls}(f_A)\in H^1(G,\af^*, {\bf PicS}(R))$$ is a   group homomorphism.\end{prop}

\begin{remark} In the proof of Proposition \ref{fi5} the following  was obtained. Let $A_1,A_2$ be two $R^\af$-Azumaya algebras containing $R$ as a maximal commutative subalgebra, and let $M_1^{(g)}, M_2^{(g)}$ be two $R$-modules such that there are isomorphism $_g(1_{g\m}A)$$\cong M_1^{(g)}\ot A_1,$ $_g(1_{g\m}A_2)\cong M_2^{(g)}\ot  A_2$ as $R\ot_{R^\af}A_1^{op}$-modules and  $R\ot_{R^\af}A_2^{op}$-modules respectively.
If $B=A_1\ot_{R^{\af}}A_2$ and $e$ is the separability idempotent of the ring extension $R/R^\af,$ then $eBe$ is an $R^\af$-Azumaya algebra containing $R$ as a maximal commutative subalgebra and  there is an $R\ot_{R^\af}(eBe)^{op}$-module isomorphism
\begin{equation}\label{ebe}_g(1_{g\m}eBe)\cong (M_1^{(g)}\ot M_2^{(g)})\ot eBe.\end{equation} Moreover, $[B] = [eBe]\in  B(R/R^\af).$

\end{remark}
By  \cite[Theorem 3.8]{DPP1}  we have  that $[M^{(g)}] \in {\bf Pic}(D_g),$ for all $g\in G.$
For further reference, we shall construct  $M^{(g)}$ explicitly.

We have the following.
\begin{lema}\label{kan}
Let A be an Azumaya $R^\af$-algebra containing R as a maximal commutative subalgebra. Write
$$C_g=C_{_{g}(1_{g\m} A)}(R)=\{c\in 1_{g\m} A  \mid \af_{g\m}( r1_{g})c =cr, \,\, r\in R\},$$
for each $g\in G.$ Then there is a left $R\ot_{R^\af}A^{op}$-module isomorphism  $_{g}(1_{g\m} A)\cong C_g\ot_RA.$

\end{lema}
\p
 Since $R/R^\af$ is a partial Galois extension then by \cite[Theorem 4.2]{DFP} the extension  $R/R^\alpha$ is separable, moreover since  $A$ contains $R$ as a maximal commutative subalgebra  we  obtain from \cite[Theorem 7.4.2]{TF} that  $A$ is a faithfully projective $R$-module and  by \cite[Excercise 7.4.7 p. 258]{TF}  that
 the functor \begin{equation*}-\ot_{R} A\colon{\rm Mod} _{R}\to\!\! _{R\ot_{R^\af} A^{op}}{\rm Mod},\end{equation*}
determines an equivalence of categories, whose inverse is given by the functor
\begin{equation*} (-)^R\colon \!\! _{(R\ot_{R^\af} A^{op})}{\rm Mod}\to {\rm Mod}_{R},\end{equation*} where
$$N^{R}=\{x\in N\mid (r\ot_{R^\af} 1)x=(1\ot_{R^\af} r)x,\, \forall\,r\in R \},$$
for each $N$ in  $ _{(R\ot_{R^\af} A^{op})}{\rm Mod}.$
In particular, for any left  $R\ot_{R^\af} A^{op}$-module $N$ there is a left   $R\ot_{R^\af} A^{op}$-module isomorphism $N\cong N^R\ot_{R} A.$  Then taking $N=\,_{g\m}(1_{g} A),$ we have
 $N^R=C_g,$ and 
the desired isomorphism follows.\cua

\begin{remark}\label{mgtocg} It follows from  Lemma \ref{kan}  that for any $[A]\in B(R/R^{\af})$  one has that $\varphi_5([A])={\rm cls}(f),$ where $f(g)=[C_g]$ for any $g\in G.$
\end{remark}

\medskip

Now we give:
\begin{lema}\label{iso2k} Let $R/R ^\af$ be a partial Galois extension. Then for an arbitrary $[J]\in {\bf Pic}(R)$  there is an $R\ot_{R^\af}    ({\rm End}_{R^\af}  (J))^{op}$-module isomorphism
 $$_{g}(1_{g\m}{\rm End}_{R^\af}(J))\cong \,_{g}(1_{g\m}J^*)\ot  J\ot( {\rm End}_{R^\af}(J)),$$ where $J^*={\rm Hom}_R(J,R),$ and the $R$-module structure of $_{g}(1_{g\m} J^*)$ is given    by (\ref{bulletaction}). 

\end{lema}
\p Let $g\in G,$ since $[J^*]$ is the inverse of $[J]$ in ${\bf Pic}(R) \subseteq {\bf PicS}(R)$ and $[D_{g\m}]$ is an idempotent in ${\bf PicS}(R)$  there are isomorphisms of $R$-modules
$$1_{g\m}J^*\cong D_{g\m}\ot J^* \cong D_{g\m}^{*}\ot J^* \cong (D_{g\m}\ot J)^*  \cong (1_{g\m} J)^*,$$
and consequently, there is an $R$-module isomorphism $_{g}(1_{g\m} J^*)\stackrel{(\ref{esg})}{\cong}( _{g}(1_{g\m}J))^*,$  where
\begin{equation}\label{iij}( _{g}(1_{g\m} J))^*=\{\eta\colon 1_{g\m} J\to R\,\mid\,\eta(\af_{g\m}(r1_{g})x)=r\eta(x),\,\, \forall r\in R, x\in   1_{g\m}J \}.\end{equation}
Now set
$C_g=\{\xi\in
1_{g} {\rm End}_{R^\af}(J)\,|\, \af_{g}(r1_{g\m})\xi(p)=\xi(rp),\,\forall p \in J, r\in R\}.$

Since $ {\rm End}_{R^\af}(J) $ is an Azumaya $R^\af$-algebra containing $R$ as a maximal commutative $R^\af$-algebra, then
by Lemma \ref{kan}  it is enough to show that there is an $R$-module isomorphism $$(_g(1_{g\m} J))^* \ot  J \cong C_g.$$

For this purpose notice that   for each $\eta\in( _{g}(1_{g\m} J))^*$  and $x\in J$   the map $$\eta_x\colon J\ni p\mt \eta(1_{g\m} p)x\in J,$$ belongs to ${\rm End}_{R^\af}(J).$
For, given $r\in R^\af$ and  $p\in J$ we have that
$$\eta_x(rp)=\eta(r 1_{g\m} p)x=\eta(\af_{g\m}(r1_{g})p)x =r\eta(1_{g\m}p)x=r\eta_x(p).$$  Furthermore, $1_{g}\eta_x\in C_g,$ because,
\begin{align*}
\af_{g}(r1_{g\m})(1_g\eta_x(p))&=\af_{g}(r1_{g\m})\eta(1_{g\m}p)x\\
&\stackrel{\eqref{iij}}= \eta[\af_{g\m}(\af_{g}(r1_{g\m})) p]x\\
&=\eta(1_{g\m}rp)x\\
&=(1_g\eta_x)(rp),\end{align*}
 for each $r\in R$ and $p\in J.$ Since $\eta_{rx}=r\eta_x,$ for all $r\in R$ and $x\in J,$  we have an $R$-linear map
$\vartheta\colon( _{g}(1_{g\m} J))^* \otimes J\ni \eta\ot x \mapsto1_{g}\eta_x\in C_g ,$ and we shall prove that $\vartheta$ is bijective.

\vu

{\bf $\vartheta$ is onto.} Consider  a dual basis $\{x_i,\eta_i\}$  for the projective $R$-module $_{g}(1_{g\m}J),$
take an arbitrary $ \xi\in C_g.$  Notice first that for any  $\eta \in (_{g}(1_{g\m} J))^*$ the map
 $$\tilde{\eta} \colon 1_{g\m}  J\ni   x\mt \af_{g\m}(1_{g}\eta(1_{g\m}\xi(x)))\in R$$
 is an element of $( _{g}(1_{g\m} J))^*$ because for any $r\in R$ and $x\in 1_gJ$ we have
 \begin{align*}\tilde\eta(\af_{g\m}(r1_{g}) x))&=\af_{g\m}\{1_{g}\eta[1_{g\m}\xi(\af_{g\m}(r1_{g})x)]\}\\
&\stackrel{\xi\in C_g} =\af_{g\m}[\eta(r1_{g\m}\xi(x))]\\
&\stackrel{\eqref{iij}}=\af_{g\m}(\af_{g}(r1_{g\m})\eta(1_{g\m}\xi(x)))\\
&=r\af_{g\m}(1_{g}\eta(1_{g\m}\xi(x)))
\\
&= r \tilde\eta( x).\end{align*}

\noindent Taking specifically the  $ \tilde\eta _i$ and  $x\in 1_gJ$  we have
 $$\tilde\eta_{i_ {x_i}}(x)=\tilde\eta_i(x)x_i=\af_{g\m}(1_{g}\eta_i(\xi(x)))x_i=
\eta _i(\xi(x))\bullet x_i.$$
From this, using the fact that   $\{x_i,\eta_i\}$ is a dual basis, we conclude that $\sum\limits_i\tilde\eta_{i_ {x_i}}(x)= \sum\limits_i\eta_i(\xi(x))\bullet x_i=\xi(x).$
Thus,
$\vartheta\left(\sum\limits_i x_i\ot \tilde\eta_i \right)=1_{g}\sum\limits_i (\tilde\eta_{i})_ {x_i}=1_g\xi=\xi,$ showing that  $\vartheta$ is surjective.

\smallskip

{\bf $\vartheta$ is injective.}
Let $c=\sum\limits_{j}\upsilon_j\ot z_j \in \ker \vartheta .$ Then $0=1_{g}\sum\limits_j (\upsilon_{j})_{z_j}.$ Consequently, for any $p\in J$ we have
 $$0=1_{g}\sum\limits_j (\upsilon_{j})_{z_j}(p) = \sum\limits_j 1_{g} \upsilon_{j}(1_{g\m}p)z_j=\sum\limits_j\upsilon_{j}(1_{g\m}p)z_j=\sum_j (\upsilon_{j})_{z_j}(p),$$ and $\sum\limits_j (\upsilon_{j})_{z_j} =0.$
  Let $\{p_l,\xi_l\}$ be a dual $R$-basis for  $J.$ Then $$c=\sum\limits_{j}z_j\ot \upsilon_j=\sum\limits_{j,l}\xi_l(z_j)p_l\ot  \upsilon_j=\sum\limits_lp_l\ot \sum\limits_j\xi_l(z_j) \upsilon_j.$$
 Since for any $x\in 1_{g\m}  J$ we have
$$\left(\sum_j\xi_l(z_j) \upsilon_j\right)(x)
=\xi_l\left(\sum_jz_j \upsilon_j(  x)\right)= \xi_l\left(\sum_j (\upsilon_j)_{z_j}(x)\right)= 0,
$$
and we conclude that $c=0$ thus $\vartheta$ is injective.\cua

 \begin{lema}\label{AigualPCP}
Let A be an Azumaya $R^\af$-algebra containing R as a maximal commutative subalgebra. If there exists a  $R\ot_{R^\alpha}A^{op}$-module isomorphism $_g(1_{g\m}A)\simeq 1_gA,$ then $A=R\star_{\alpha,\omega}G$, for some $\omega \in H^2(G,\alpha,R)$.
\end{lema}

\p Let $\0_g\colon _g(1_{g\m} A) \to 1_g  A$ be an $R\ot_{R^\af} A^{op}$-module isomorphism. Write $\0_g(1_{g\m})=u_g.$ 
We shall construct an $R^\af$-algebra $\B$ which is isomorphic to a crossed product by a twisting of $\af$. To this  we establish first some properties of the family $\{u_g\}_{g\in G}.$ For  $r\ot_{R^\af} a\in R\ot_{R^\af}A^{op},$  we have
$$(r\ot_{R^\af} a)\bullet 1_{g\m}=\af_{g\m}(r1_g)a=(1\ot_{R^\af} \af_{g\m}(r1_g)a)\bullet (1_{g\m}).$$
Applying $\0_g$ we get
$$(r\ot_{R^\af} a) u_g=(1_g\ot_{R^\af} \af_{g\m}(r1_g)a)  u_g,$$ or,
$r u_ga=u_g\af_{g\m}(r1_g)a$ In particular, letting $a=1$  we get
\begin{equation}\label{ug} r u_g=u_g\af_{g\m}(r1_g),\,\,\, \text{for any}\,\, \,r\in R\,\, \,\text{and}\,\,\, g\in G.
\end{equation}
It follows from (\ref{ug}) that
\begin{equation}\label{uu}
u_g=u_g1_{g\m}\,\,\,\,\text{ for all} \,\,\,\,g\in G.\end{equation}
Now, take $a_g\in 1_{g\m} A,$  such that $\0_g(a_g)=1_g.$ Then,
$1_g= \0_g((1_g\ot a_g)\bullet (1_{g\m}))=(1_g\ot a_g)  \0_g(1_{g\m})
=   u_ga_g,$ and we conclude that
\begin{equation}\label{invd}1_g=u_ga_g,\,\,\,\,\,\text{for all}\,\,\,\, g\in G.
\end{equation} Moreover,   $\0_g( a_gu_g)=(1_g\ot u_g)\0_g( a_g)=u_g,$ and since $\0_g$ is injective we get 
\begin{equation}\label{invi}1_{g\m}= a_gu_g,
\end{equation} and follows from \eqref{ug} that for each $r\in R$ we have  $ a_gr u_g=a_gu_g\af_{g\m}(r1_g)\stackrel{\eqref{invi}}=\af_{g\m}(r1_g),$ which gives 
$$a_gr1_g\stackrel{\eqref{invd}}=(a_gr u_g)a_g=\af_{g\m}(r1_g)a_g,$$ since $a_g\in 1_{g\m}A$ we get
\begin{equation}\label{agc}
{a}_{g}r =
\af_{g\m}(r1_{g}){a}_{g},
\end{equation}for all $g\in G.$
Let $w_{g,h}=u_gu_h{a}_{gh},$ for all $g,h\in R.$ Notice that $w_{g,h}\in R,$ because for each $r\in R$
\begin{align*}rw_{g,h}&=ru_gu_h{a}_{gh}
\\&\stackrel{(\ref{ug})}=u_gu_h\af_{h\m}(\af_{g\m}(r1_g)1_h){a}_{gh}\\
&=u_gu_h\af_{(gh)\m}(r1_{gh})1_{h\m}{a}_{gh}\\
&\stackrel{(\ref{uu})}{=}u_gu_h\af_{(gh)\m}(r1_{gh}){a}_{gh}
\\&\stackrel{(\ref{agc})}{=}u_gu_h{a}_{gh}r.\end{align*}

Now we check that $w_{g,h}\in D_{g}D_{gh},$ indeed
\begin{align*}
w_{g,h}1_g1_{gh}&=1_gu_gu_h{a}_{gh}1_{gh}
\\&=u_gu_h{a}_{gh}1_{gh}
\\&\stackrel{(\ref{agc})}=u_gu_h1_{(gh)\m}{a}_{gh}
\\&\stackrel{(\ref{uu})}=u_gu_h1_{h\m}1_{(gh)\m}{a}_{gh}
\\&=u_gu_h\af_{h\m}(1_{g\m}1_h){a}_{gh}
\\&\stackrel{\eqref{ug}}=(u_g1_{g\m})u_h{a}_{gh}
\\&\stackrel{\eqref{uu}}=u_gu_ha_{gh}.
\end{align*}
In an analogous way one can show that $r_{g,h}=u_{gh}a_ha_g\in D_gD_{gh}.$
Thus, for any $g,h\in G$ we see that
\begin{align*}
w_{g,h}r_{g,h}&=u_gu_h{a}_{gh}u_{gh}a_ha_g
\\&=u_gu_h1_{(gh)\m}a_ha_g
\\&=u_gu_ha_ha_g
\\&=u_g1_ha_g
\\&=u_g1_h1_{g\m}a_g
\\&=u_g\af_{g\m}(1_{gh}1_{g})a_g
\\&=1_{gh}1_{g}u_ga_g.
\\&=1_{gh}1_{g}.
\end{align*}

Hence to prove that the family ${w}=\{w_{g,h}\}_{(g,h)\in G\times G}$ is a twisting for $\af$ we need to check  equality (v).
For any $g,h,l\in G$ we have
$$ w_{g,h}w_{gh,l}=u_gu_h{a}_{gh}u_{gh}u_l{a}_{ghl}
\stackrel{(\ref{invi})}=u_gu_h1_{(gh)\m}u_l{a}_{ghl}=u_gu_hu_l{a}_{ghl},
$$
and
$$\af_g(w_{h,l}1_{g\m})w_{g,hl}=\af_g(r_{h,l}1_{g\m})u_gu_{hl}{a}_{ghl}
\stackrel{(\ref{ug})} =u_gu_hu_l{a}_{hl}u_{hl}{a}_{ghl}\stackrel{(\ref{invi})} =u_gu_hu_l{a}_{ghl}.$$
\smallskip

Set $\B=\sum\limits_{g\in G}\B_g\subseteq A,$  where $\B_g=D_gu_g,$ for any $g\in G.$ Notice that by \eqref{ug} we get $\af_{g}(r1_{g\m}) u_g=u_gr,$ for all $r\in R$ and  $g\in G,$ moreover
$w_{g,h}u_{g,h}=u_gu_h1_{(gh)\m}=u_gu_h.$ Then by \eqref{product} the map $$\digamma\colon R\star_{\af,w} G\ni\sum\limits_{g\in G} a_g\de_g\to \sum\limits_{g\in G} a_gu_g \in \B$$
is an $R^\af$-algebra epimorphism. Since $R\star_{\af,w} G$ is an Azumaya $R^\af$-algebra, \cite[Corollary 2.3.7]{DI}
implies   $\ker \digamma=(\ker \digamma\cap R^\af)(R\star_{\af,w} G).$ Moreover, the fact that $u_1 \in {\mathcal U}(R),$   yields that   $\digamma$ restricted to $R^\af$ is injective, and, consequently,   $\digamma$ is a ring monomorphism. This yields that     $\B=\bigoplus\limits_{g\in G}\B_g\cong R\star_{\af,\mathfrak{r}} G$ as $R^\af$-algebras, and we may assume that  $R\star_{\af,w} G\subseteq A.$

\vu

We  check that $C_A(R\star_{\af,w}G)=R^\af.$ From $R\subseteq R\star_{\af,\mathfrak{r}}G, $ we obtain that
$$R^\af\subseteq C_A(R\star_{\af,w}G)\subseteq C_A(R)=R.$$
Take $x\in C_A(R\star_{\af,w}G),$ then $x\in R$ and for all $g\in G$ we have $x(1_g\de_g)=(1_g\de_g)x,$ or equivalently $x1_g=\af_g(x1_{g\m}),$ which gives $x\in R^\af,$ showing that $C_A(R\star_{\af,w}G)=R^\af.$

\vu

Since $R\star_{\af,w}G$ is an Azumaya $R^\af$-subalgebra of $A,$   the double centralizer Theorem (see \cite[Theorem 2.4.3]{DI}) implies  $R\star_{\af,w}G=C_A(C_A(R\star_{\af,w}G))=C_A(R^\af)=A.$\cua

\begin{teo}\label{fi4fi5} The sequence $H^2(G,\af, R)\stackrel{\varphi_4}{\to} B(R/R^\af)\stackrel{\varphi_5}{\to}  H^1 (G,\af^*,{\bf PicS}(R)) $ is exact.
\end{teo}

\p Let ${\rm cls}(\om)\in H^2(G,\af, R)$ and $\varphi_5([R\star_{\af,\om}G])={\rm cls}(f).$ We shall check that $f\in B^1(G,\af^*, {\bf PicS}(R))$ by proving that $f(g)=[D_g],$ for any $g\in G.$ For this we need to show  that $_g (D_{g\m}\ot R\star_{\af,\om}G) $ and $ D_g\ot R\star_{\af,\om} G$ are isomorphic as  $ D_g\ot  _{R^{\af }}  (R\star_{\af,\om}G)^{op}$-modules.

 Notice that the map $$D_{g\m}\times R\star_{\af,\om}G \ni\left(d, \sum_{h\in G}a_h\de_h\right) \mapsto 1_g\ot\sum_{h\in G}\afg(da_h)\om_{g,h}\de_{gh}\in D_g\ot R\star_{\af,\om}G $$ is $R$-balanced, so that it induces a well defined map

$$\nu_g\colon _g(D_{g\m}\ot R\star_{\af,\om}G) \ni d\ot \sum_{h\in G}a_h\de_h \mapsto 1_g\ot \sum_{h\in G}\afg(da_h)\om_{g,h}\de_{gh}\in D_g\ot R\star_{\af,\om}G.$$

Now for $d_g\ot  _{R^{\af }} a_l\de_l\in  D_g\ot _{R^{\af }} (R\star_{\af,\om}G)^{op}$ and  $d\ot \sum_{h\in G}a_h\de_h\in \,_g(D_{g\m}\ot R\star_{\af,\om}G),$ we have
\begin{align*}&(d_g\ot  a_l\de_l)\bullet \left(d\ot \sum_{h\in G}a_h\de_h\right)=\af_{g\m}(d_g)d\ot \sum_{h\in G}a_h\afh(a_l1_{h\m})\om_{h,l}\de_{hl}\\
&\stackrel{\nu_g}{\mapsto}1_g\ot \sum_{h\in G}d_g\afg(da_h\afh(a_l1_{h\m})\om_{h,l})\om_{g,hl}\de_{ghl}\\
&\stackrel{(\ref{prodp})}{=}d_g\ot \sum_{h\in G}\afg(da_h)\af_{gh}(a_l1_{(gh)\m})\afg(\om_{h,l}1_{g\m})\om_{g,hl}\de_{ghl}\\
&=d_g\ot\sum_{h\in G}\afg(da_h)\om_{g,h}\af_{gh}(a_l1_{(gh)\m})\om_{gh,l}\de_{ghl}\\
&=d_g\left(1_g\ot\sum_{h\in G}\afg(da_h)\om_{g,h}\de_{gh} \right)(a_l\de_l)\\
&=(d_g\ot a_l\de_l) { \nu _g } \left(d\ot \sum_{h\in G}a_h\de_h\right),
\end{align*}
and $\nu_g$ is $D_g\ot _{ R^{\af} } (R\star_{\af,\om}G)^{op}$-linear.  To prove that $\nu_g$ is injective consider the map
$$ D_g\ot R\star_{\af,\om}G \ni d' \ot \sum_{h\in G}a_h\de_h \stackrel{\lb_g}\mt 1_{g\m}\ot \sum_{h\in G}\af_{g\m}(d'a_h \, \om\m_{g,h} ) \de_{g\m h}\in\,  _g(D_{g\m}\ot R\star_{\af,\om}G).$$
It is easy to see that $\lb_g\nu_g={\rm id}_{_{g}(D_{g\m}\ot R\star_{\af,\om}G)},$ and $\nu_g$ is injective.

\vu

For the surjectivity, consider $d_g\ot\sum_{h\in G}a_{gh}\de_{gh}=1_g\ot\sum_{h\in G}d_ga_{gh}\de_{gh}\in D_g\ot R\star_{\af,\om}G.$ Then,
$$\nu_g\left(\af_{g\m}(d_g)\ot\sum_{h\in G}\af_{g\m}(a_{gh}\om\m_{g,h})\de_h\right)=1_g\ot \sum_{h\in G}d_ga_{gh}\om\m_{g,h}\om_{g,h}\de_{gh}=d_g\ot\sum_{h\in G}a_{gh}\de_{gh}.$$
and $\nu_g$ is surjective. Using (\ref{unico}) we get   $f(g)=[D_g],$ which implies $\varphi([R\star_{\af,\om}G])={\rm cls}(1)$ in $H^1(G,\af^*,{\bf PicS}(R))$ and ${\rm im}\,\varphi_4\subseteq \ker\varphi_5.$
\medskip

Now we prove ${\rm im}\,\varphi_4\supseteq \ker\varphi_5.$
Let $[A]\in \ker\varphi_5$ with $R$ being a maximal commutative subalgebra of $A.$ Then, $f=f_A\in B^1(G,\af^*, {\bf Pic}S(R)),$ where $f_A$ is given by \eqref{fa},  and there exists $[P]\in {\bf Pic}(R)$ such that  $M^{(g)}\cong (1_{g\m} P)_{g}\ot  P^*$ as $R$-modules. By the construction of $\f_5,$ the latter implies that there is  an $R\ot_{R^\af} A^{op}$-module isomorphism
\begin{equation}\label{ISO}\,_g(1_{g\m} A)\cong _g(1_{g\m} P)\ot P^*\ot A.\end{equation}
Since $D_g$ is a direct summand of $R,$ the map $  {\bf Pic}(R)\ni [Q] \to [D_g\ot Q]\in {\bf Pic}(D_g)$  is a group epimorphism, for each $g\in G$. Consequently,    $[1_{g\m}P]=[D_{g\m}\ot P]\in {\bf Pic}(D_{g\m})$ and $  [_g( 1_{g\m} P) ] \in {\bf Pic}(D_{g})$   by  the first item of  Lemma \ref{gaction}. Using again the above epimorphism,  we see
 that there exists $[J]\in {\bf Pic}(R)$ and a $D_g$-module isomorphism $ _{g}( 1_{g\m}  P) \cong D_{g}\ot  J^\ast\cong 1_gJ^*.$

 Applying ${\af_g }^*$ to both parts of the equality we obtain an $R$-module isomorphism   $1_{g\m}P\cong  _{g\m}(1_{g} J^*),$ Moreover,   $(1_{g\m} P^*)_{g}\cong 1_gJ,$ as $R$-modules, in view of (\ref{esg}).   By Lemma \ref{iso2k}  there exists a $R\ot_{R^\af} ({\rm End}_{R^\af}(J))^{op}$-module isomorphism
\begin{equation*}\label{isook}_g(1_{g\m}{\rm End}_{R^\af}(J))\cong \,_{g\m}(1_g J^*)\ot  J\ot {\rm End}_{R^\af}(J).
\end{equation*}
 Since $_g(1_{g\m}  {\rm End}_{R^\af}(J))$ is a left unital $D_g$ module, tensoring by $D_g$ we obtain
\begin{align*}_g(1_{g\m} {\rm End}_{R^\af}(J)) & \cong D_g \ot \,_{g\m}(1_{g}  J^*)\ot   J\ot {\rm End}_{R^\af}(J)\\
& \cong D_g \ot 1_{g\m}P\ot   J\ot {\rm End}_{R^\af}(J)\\
& \cong\, \,1_{g\m}P\ot  D_g \ot  J\ot  {\rm End}_{R^\af}(J)\\
& \cong\, \,1_{g\m}P\ot  1_g  J\ot  {\rm End}_{R^\af}(J)
\\&\cong  1_{g\m}P\ot \, _g(1_{g\m} P^*)\ot {\rm End}_{R^\af}(J),\end{align*}\\
as $R\ot_{1_gR^\af} {\rm End}_{R^\af}(J)^{op}$-modules.
Let $e$ be an idempotent of separability of $R$ over $R^\af.$ It follows from  \eqref{ebe}, \eqref{ISO}  and the latter isomorphism, that taking $\Delta=A \ot_{R^\af} {\rm End}_{R^\af}(J) $ one has as $R\ot_{R^\af}\Delta^{op}$-modules that
\begin{align*}  _g( 1_{g\m }  e\Delta e )
\cong&
 \{ _g(1_{g\m}  P)\ot P^* \}   \ot \{ 1_{g\m} P\ot  _g(1_{g\m}P^*) \} \ot    e\Delta e \\
\cong&
 \{ _g(1_{g\m}  P)\ot ( P^*    \ot  1_{g\m} P)\ot  _g(1_{g\m}P^*) \} \ot   e\Delta e
\\
\cong&
 \{ _g(1_{g\m}  P)\ot D_{g\m}\ot  _g(1_{g\m}P^*) \} \ot  e\Delta e \\
\cong&
 \{ _g(1_{g\m}  P)\ot_g(1_{g\m}P^*) \} \ot   e\Delta e \\
 \cong & _{g}(D_{g\m})\ot     e\Delta e
\\ \cong &D_{g}\ot    e\Delta e.
\end{align*}
Since  $[e\Delta e] =  [\Delta ] =[A]$  in $B(R^\af),$ replacing $e\Delta e$ by $A,$ we may suppose,  without loss of generality, that
 $_g(1_{g\m}A)  \cong 1_g  A, $
as $R\ot_{R^\af}A^{op}$-modules. By Lemma  \ref{AigualPCP},  we have that  $A=R\star_{\alpha,\omega}G$, for some  $\omega \in H^2(G,\alpha,R)$. Then, $\varphi_4(cls(\omega))=[R\star_{\alpha,\omega} G]=[A]$, and  $[A]\in \mbox{im}\varphi_4.$\cua

\section{The homomorphism $\varphi_6\colon H^1(G,\af^*, {\bf PicS}(R))\to H^3 (G,\af, R)  $}  We start by recalling    from \cite{DPP1} the definition of the   $R$-$R$-bimodule $ _{g}{(D_{g\m})}_{I},$  $g\in G:$   its  underlying set is $D_{g\m},$ endowed with the   action $*$ given by
$$r*d=\af_{g\m}(r1_{g})d,\,\,\,\,\text{and} \,\,\, d*r=dr, \,\,\, \text{for any}\,\,\, r\in R,\,d\in D_{g\m}.$$
Write for simplicity $_{g}{(D_{g\m})}=\, _{g}{(D_{g\m})}_{I}.$
Now we give from \cite{DEP} the concept of a partial representation.

\begin{defi}\label{defn-par-repr} A (unital) partial representation of $G$ into an algebra (or, more ge\-ne\-ral\-ly, a monoid) $S$ is a map $\Phi: G \to S$ which satisfies the following properties, for all $g,h\in G,$

\vspace*{2mm}
\noindent (i) $\Phi (g\m) \Phi (g) \Phi (h) = \Phi (g\m) \Phi (g h),$

\vspace*{2mm}
\noindent (ii) $\Phi (g ) \Phi (h) \Phi (h\m) = \Phi (g h) \Phi ( h \m ),$

\vspace*{2mm}
\noindent (iii) $\Phi (1_G )  = 1_S.$
\end{defi}

 We also recall  two important for us partial representations of $G$.
\begin{prop}\label{pr1} Let $\Phi_0\colon G\ni g\mt [_g(D_{g\m})]\in {\bf PicS}_{R^\af}(R)$  and for  any cocycle $f\in Z^1(G,\af^*, {\bf PicS}(R)),$  set $\Phi_f=f\Phi_0\colon G\to  {\bf PicS}_{R^\af}(R),$ that is $\Phi_f(g)=f(g)\Phi_0(g),$ for any $g\in G.$   Then,
\begin{itemize}
\item \cite[Proposition 6.2]{DPP1}  $\Phi_0$  is a partial representation   with $\Phi_0(g)\Phi_0(g\m)=[D_g].$
\item \cite[Lemma 6.3]{DPP1}
 $\Phi_f$ is a partial representation of G, with  $\Phi_f(g)\Phi_f(g\m)=[D_g].$
 Moreover,
 writing
$\Phi_f(g)=[J_g],$  we have that  $D_g\cong {\rm End}_{D_g}(J_g),$ as $R$-
and   $D_g$-algebras,  for any $g\in G.$

\end{itemize}
\end{prop}

\begin{remark}\label{jgdg} Let $f$ be an element of $ Z^1(G,\af^*, {\bf PicS}(R))$ and  $J_g=f(g)\ot _g(D_{g\m}).$ Then
\begin{itemize} \item \cite[Remark 6.5]{DPP1} For  $ x_g\in J_g$ and  $r\in R$ we have.
\begin{equation}\label{jgtwits}\af_g(r1_{g\m})x_g=x_gr,\,\, \end{equation} and if  $D_g$ is a semi-local ring,  then for any $g\in G,$ we see   that there is   $u_g\in J_g,$ a $D_g$-basis of $J_g$ with  $J_g=D_gu_g,$ and  $u_gr=\af_g(r1_{g\m})u_g,$ for all $ r\in R,\,g\in G.$

\item By \cite[Lemma 6.7]{DPP1}  the map $ J_g\ot  D_h \stackrel{\kappa_{g,h}}\to D_{gh}\ot  J_g$ such that
$$a_{g}\ot \, b_h\mt\af_g(b_h1_{g\m})\ot a_{g},$$
for any $g,h\in G,$ is a $R$-$R$-bimodule isomorphism.
\end{itemize}

\end{remark}

{\bf The construction of $\varphi_6.$}
Take  $f\in Z^1(G,\af^*, {\bf PicS}(R)).$ By the second item of Proposition \ref{pr1}
 there exists a family of $R$-$R$-bimodule isomorphisms $\{f_{g,h}\colon J_g\ot J_h\to D_g\ot J_{gh}\}_{g,h\in G}.$
 Consider the following diagram
\begin{equation} \label{diag1}
\xymatrix{J_g\ot J_h\ot  J_l\ar[d]^{ f_{g,h}\ot{\rm id}_l}\ar[r]^{{\rm id}_g\ot f_{h,l}}& J_g\ot  D_{h}\ot J_{hl} \ar[r]^{\kappa_{g,h}\ot {\rm id }_{hl}}&\,\,\,D_{gh}\ot J_g \ot J_{hl}\,\,\ar[r]^{{\rm id}_{D_{gh}}\ot f_{g,hl}} &\,\,\,\,D_{gh}\ot  D_g\ot  J_{ghl}\ar[d]^{{\tau}_{gh,g}\ot {\rm id}_{ghl}} \\
D_g\ot J_{gh}\ot J_l\ar[rrr]^{{\rm id}_g\ot f_{gh,l}}&&&D_{g}\ot D_{gh}\ot J_{ghl}},
\end{equation}
for any $g,h,l\in G,$ where $\kappa_{g,h}$ is given by the second item of  Remark \ref{jgdg}.
Set $\tilde\om(g,h,l)=({\rm id}_g\ot f_{gh,l})\circ(f_{g,h}\ot{\rm id}_l)\circ({\rm id}_g\ot f_{h,l})\m\circ(\kappa_{g,h}\ot {\rm id }_{hl})\m\circ({\rm id}_{D_{gh}}\ot f_{g,hl})\m\circ({\tau}_{gh,g}\ot {\rm id}_{ghl})\m.$

\vu

Then,  $$\tilde\om(g,h,l)\in \U({\rm End}_{D_g\ot  D_{gh}\ot  D_{ghl}}(D_g\ot D_{gh}\ot  J_{ghl})),$$
 and there is a unique $\om_1(g,h,l)\in \U(D_g\ot D_{gh}\ot D_{ghl})$ such that
 $\tilde\om(g,h,l)(z)=\om_1(g,h,l)z$\,\, \text {for all}\,\, $z\in D_g\ot D_{gh}\ot J_{ghl}.$  This implies  that there is a unique  element $\om_f(g,h,l)\in \U(D_gD_{gh}D_{ghl})$ with $$\tilde\om(g,h,l)z=\om_f(g,h,l)z, g,h,l\in G,\, z\in D_g\ot  D_{gh}\ot J_{ghl}. $$

\begin{prop} \cite[Theorem 6.9]{DPP1}\label{isof6}  The map $H^1(G,\af^*, {\bf PicS}(R))\ni {\rm cls}(f)\stackrel{\f_6}{\to} {\rm cls}(\om_f)\in  H^3(G,\af, R) $ is a homomorphism of groups.
\end{prop}

\vu

\section{Partial generalized crossed products and exactnesss at\\  $H^1(G,\af^*, {\bf PicS}(R))$ }\label{pargencrosprod}

Let $f\in Z^1(G,\af^*, {\bf PicS}(R))$ and  $\Phi_f=f\Phi_0$. Put $\Phi_f(g)=[J_g], \, g\in G,$ and suppose that the diagram (\ref{diag1}) given by the family of $R$-$R$-bimodule isomorphisms $\mathfrak{F}=\{f_{g,h}\colon J_g\ot J_h\to D_g\ot J_{gh}\}_{g,h\in G}$ commutes. Following \cite{K} we call $\mathfrak{F}$ a {\it factor set of G related to  $\Phi = \Phi _f.$}  Consider the $R$-$R$-bimodule
$$\Delta=\Delta(\mathfrak{F},\af,R,\Phi, G)=\bigoplus\limits_{g\in G}J_g.$$ Let $m_{g,h}\colon D_g\ot J_{gh}\to D_gJ_{gh}\subseteq J_{gh},g,h\in G,$ be the multiplication map\footnote{Notice that if $g=1$ or $h=1$ then $m_{g,h}$ is an $R$-$R$-bimodule isomorphism.}. The product of elements in $\Delta$ is defined by the formula
$xy=m_{g,h}\circ f_{g,h}(x\ot y),\,\, \,\text{for any}\,\, \,x\in J_g, \,\, y\in J_h.$

We call $\Delta(\mathfrak{F},\af,R,\Phi, G)$ a {\it partial generalized crossed product.} Notice that for any $x\in J_g,y\in J_h,$ and $f_{g,h}(x\ot y)=\sum d_g\ot u_{gh}$ we have $xy=\sum d_g u_{gh}.$ From this we obtain
\begin{equation}\label{mulf}f_{g,h}(x\ot y)=1_g\ot \sum d_gu_{gh}=1_g\ot xy.\end{equation}
In particular, for any $r\in R,$ $r(xy)=\sum rd_g u_{gh}=m_{g,h}\circ f_{g,h}(rx\ot y)=(rx)y,$ and $(xy)r=m_{g,h}\circ f_{g,h}(x\ot yr)=x(yr).$ We conclude that
\begin{equation}\label{asp} r(xy)=(rx)y, \,\, {\rm and}\,\,\, (xy)r=x(yr)\,\,\,\text{for any}\,\, g,h\in G,\,\,r\in R,\,\, x\in J_g\,\, \text{and}\,\, y\in J_h.
\end{equation}

\vu

\begin{prop} \label{pgcp1}The partial generalized crossed product $\Delta=\Delta(\mathfrak{F},\af,R,\Phi, G)$ is an associative $R^\af$-algebra with identity element and $J_1\cong R$ as $R$-algebras.
\end{prop}
\p Let $x\in J_g,y\in J_h, z\in J_l.$ Then
\begin{align*}
1_g\ot 1_{gh}\ot x(yz)&=(\tau_{gh,g}\ot \iota_{ghl})(1_{gh}\ot 1_g\ot  x(yz))\\
&\stackrel{(\ref{mulf})}{=}(\tau_{gh,g}\ot \iota_{ghl})(1_{gh}\ot f_{g,hl}(x\ot yz))\\
&=(\tau_{gh,g}\ot \iota_{ghl})({\rm id}_{D_{gh}}\ot  f_{g,hl})(1_{gh}\ot x\ot yz)\\
&\stackrel{x\in J_g}=(\tau_{gh,g}\ot \iota_{ghl})({\rm id}_{D_{gh}}\ot  f_{g,hl})(1_g1_{gh}\ot x\ot yz)\\
&=(\tau_{gh,g}\ot \iota_{ghl})({\rm id}_{D_{gh}}\ot  f_{g,hl})(\kappa_{g,h}\ot \iota_{hl})(x\ot 1_h\ot yz)\\
&\stackrel{(\ref{mulf})}{=}(\tau_{gh,g}\ot \iota_{ghl})({\rm id}_{D_{gh}}\ot  f_{g,hl})(\kappa_{g,h}\ot \iota_{hl})(\iota_g\ot f_{h,l})(x\ot y\ot z)\\
&\stackrel{(\ref{diag1})}{=}({\rm id}_{D_g}\ot f_{gh,l})(f_{g,h}\ot \iota_l)(x\ot y\ot z)\\
&=({\rm id}_{D_g}\ot f_{gh,l})(f_{g,h}(x\ot y)\ot z)\\
&\stackrel{(\ref{mulf})}{=}({\rm id}_{D_g}\ot f_{gh,l})(1_g \ot xy\ot z)\\
&=1_g \ot f_{gh,l}(xy\ot z)\\
&\stackrel{(\ref{mulf})}{=} 1_g\ot 1_{gh}\ot (xy)z,
\end{align*}
for any $g,h\in G,$ which implies $(1_g1_{gh})[x(yz)]=(1_g1_{gh})[(xy)z].$
Furthermore, since $x\in J_g$ and $y\in J_h,$ using (\ref{asp}) we have

\begin{align*}
(1_g1_{gh})[x(yz)]&=(1_{gh}x)(yz)
\\&\stackrel{(\ref{jgtwits})}=(x1_h)(yz)
\\&=(m_{g,hl}\circ f_{g,hl})(x1_h\ot yz)
\\&=m_{g,hl}\circ f_{g,hl}(x\ot yz)
\\&= x(yz).
\end{align*}
  Moreover, since $xy\in J_{gh},$ we get
\begin{align*}
(1_g1_{gh})[(xy)z] & \stackrel{(\ref{asp})} = [(1_g1_{gh})(xy)]z
\\&= [(1_g(xy)]z\\
&\stackrel{(\ref{asp})} =  [(1_g x) y ]z
\\&= (xy)z,
\end{align*} and $\Delta$ is  associative.

\vu

Now we will check that $J_1\cong R$ as $R$-algebras.
For $x,y\in J_1, xy=m_{1,1}f_{1,1}(x\ot y)\in J_1,$ and $J_1$ is closed by products. Moreover, $[J_1]=\Phi(1)=[R],$ and $J_1\cong R$ as $R$-$R$-bimodules. Let $\phi\colon R\to J_1$ be  an  $R$-$R$-bimodule isomorphism. Denoting $u=\phi(1)$ we have $J_1=\phi(R)=Ru=uR,$ and  it follows from (\ref{jgtwits}) that  $ur=ru$, for any $r\in R$.

\vu

 Since $m_{1,1}\circ f_{1,1}$ is a $R$-$R$-bimodule isomorphism, $u\ot u$ is a basis of $J_1\ot J_1$ and then there exists $c\in\U(R)$ such that $m_{1,1}\circ f_{1,1}(u\ot u)=cu.$ If we set $e=c\m u\in J_1,$ then $J_1=eR,$ $ee=m_{1,1}\circ f_{1,1}(c\m u\ot c\m u)=e$ and the map $R\ni r\mt re\in J_1$ is a $R$-module isomorphism.

 \vu

The element $e$ is the identity of $\Delta.$ Indeed, let $x\in J_l, l\in G.$ Since $m_{1,l}\circ f_{1,l}\colon J_1\ot J_l\cong J_l$ is a $R$-$R$-bimodule isomorphism, there exists $\sum\limits_{i} r_ie\ot y_i   \in  J_1\ot J_l,$ such that $x=m_{1,l}\circ f_{1,l}(\sum\limits_{i} r_ie\ot y_i).$ Writing $y=\sum\limits_{i} r_iy_i\in J_l$ we have that $\sum\limits_{i} r_ie\ot y_i= e \ot y$ and   $x=m_{1,l}\circ f_{1,l}(e\ot y)=ey.$ Analogously using the $R$-$R$-bimodule isomorphism $m_{l ,1 }\circ f_{l,1 }\colon J_l\ot J_1 \cong J_l$ we see that  there exists $y'\in J_l$ such that $x=y'e.$ Then  $ex=e(ey)=(ee)y=ey=x$ and $xe=(y'e)e=y'e=x$.

\vu

 Therefore,  $\Delta$ is an $R^\af$-algebra with $1_{\Delta}=e,$   and the map $R\ni r\mt re\in J_1$ is multiplicative, which yields that  $J_1\cong R$ as $R$-algebras.\cua

 \vu

\begin{exe} \label{cpp} Suppose that  $f(g)=[D_g]$ for all $g\in G.$  Then   $\Delta(\mathfrak{F},\af,R,\Phi, G)$ is a partial crossed product with respect to  a twisting of $\af.$
Indeed,  as in Remark \ref{jgdg} we have $J_g=D_gu_g,$ where $u_gr=\af_g(r1_{g\m})u_g,$ and $u_g$ is a free $D_g$-basis of $J_g$  for all $g\in G.$ Write $f_{g,h}(u_g\ot u_h)=\om_{g,h}\ot u_{gh} \in D_g D_{gh} \ot u_{gh } = D_g \ot D_{gh} u_{gh}.$  Then, $ \om_{g,h}  \in \U (D_g D_{gh}).$ In fact,   since $f_{g,h}$ is an isomorphism of $R$-$R$-bimodules, there exists $x= \sum _i a_g^i u_g \ot b_h^i u_h \in D_g u_g \ot D_h u_h$ with  $ f_{g,h} (x) =   1_g1_{gh} \ot u_{gh}.$ Then, $$x=  \sum _i a_g^i u_g b_h  \ot  u_h
 \stackrel{(\ref{jgtwits})}= \sum_i a^i_g \af _g (b^i_h 1_{g\m }) u_g \ot u_h =  x_{g,h} u_g \ot u_h,$$ where $   x_{g,h}= \sum_i a^i_g \af _g (b^i_h 1_{g\m }) \in D_g D_{gh}.$ Consequently,
$$  1_g1_{gh} \ot u_{gh} =  f_{g,h} (x) =  x_{g,h}  f_{g,h}( u_g \ot u_h) = x_{g,h} \om_{g,h}\ot u_{gh} ,$$ which gives  $ x_{g,h} \om_{g,h} u_{gh} =  1_g1_{gh}  u_{gh}.$ Since $u_{gh}$ is a free generator of $J_{gh}$ over $D_{gh},$ we conclude that    $ x_{g,h} \om_{g,h}  =  1_g1_{gh} ,$ which shows that $  \om_{g,h}  \in \U (D_g D_{gh}).$

\vu

Notice that
 $$u_gu_h=m_{g,h}\circ f_{g,h}(u_g\ot u_h)=\om_{g,h}u_{gh}.$$
 Since $(u_gu_h)u_l=u_g(u_hu_l),$ for any $g,h,l\in G,$ it follows that  the map $\om\colon G\times G\ni (g,h)\to \om_{g,h}\in R$ belongs to $Z^2(G,\af,R).$   We conclude by (\ref{jgtwits})  that  $\Delta\cong \bigoplus\limits_{g\in G}D_gu_g=R\star_{\af, \om}G$ as $R$-$R$-bimodules and $R^\af$-algebras.\end{exe}

By Proposition \ref{pgcp1} one may identify $R\cong J_1\subseteq \Delta.$ In particular, $1_R=1_\Delta.$ This will be assumed in all what follows.

Our   next result states that every partial generalized crossed product is an Azumaya $R^\af$-algebra which is split by $R.$

\vu

\begin{prop} \label{gpcpisa}Let $f \in Z^1(G,\af^*,{\bf PicS}(R)),$ $\Phi=f\Phi_0$ and $\mathfrak{F}$ a factor set of $G$ related to $\Phi.$ Then the partial generalized crossed product $\Delta=\Delta(\mathfrak{F},\af,R,\Phi, G)$ is an Azumaya $R^\af$-algebra containing $R$ as a maximal commutative subalgebra.
\end{prop}
\p For $h\in G,$ we see that
\begin{equation}\label{carjh}
J_h=\{z\in \Delta\mid zr=\af_{h}(r1_{h\m})z,\,\,\, \text{for all}\,\,\, r\in R\}.
\end{equation}   Indeed,  fix $h\in G.$ Then, for any   $z=\sum\limits_{g\in G}z_g\in \Delta ,$ satisfying  $zr=\af_{h}(r1_{h\m})z$ for all $r\in R,$ we have
$$\sum\limits_{g\in G}\af_{h}(r1_{h\m})z_g=\af_{h}(r1_{h\m})z=zr=\sum\limits_{g\in G}z_gr=\sum\limits_{g\in G}\af_g(r1_{g\m})z_g,$$
 and
\begin{equation}\label{afgtoh}\af_{h}(r1_{h\m})z_g=\af_g(r1_{g\m})z_g,\end{equation} for all $r\in R,$  $g\in G.$ In particular, for $r_h=\af_{h\m}(r1_h)$ we  obtain $$(r1_h)z_g=\af_h(r_h)z_g=\af_g(r_h1_{g\m})z_g\stackrel{(\ref{prodp})}=\af_{gh\m}(r1_{h g\m})1_gz_g=\af_{gh\m}(r1_{h g\m})z_g,$$
for any $g\in G$ and any $r\in R.$ Since $R/R^\af$ is a partial Galois extension, there exists a family $\{x_i, y_i\,|\, 1\le i\le m\}\subseteq R,$   such that $\displaystyle\sum_{i=1}^{m}x_i\af_{gh\m}(y_i 1_{hg\m})=\delta_{1,gh\m},$ for any $g\in G.$  \\
Thus for $g\ne h,$
$$1_hz_g=\displaystyle\sum_{i=1}^{m}x_iy_i (1_hz_g)=\displaystyle\sum_{i=1}^{m}x_i[(y_i 1_h)z_g]=\displaystyle\sum_{i=1}^{m}x_i\af_{gh\m}(y_i 1_{hg\m})z_g=0,$$ and $1_hz_g=0.$ From this we get
$$z_g=1_gz_g=\af_{g}(1_R1_{g\m})z_g\stackrel{\eqref{afgtoh}}=\af_h(1_R1_{h\m})z_g=1_hz_g=0.$$ Therefore $z=z_h\in J_h$ and using \eqref{jgtwits} we obtain (\ref{carjh}).
In particular, if $z\in C_\Delta(R)$ we have $zr=rz$ for any $r\in R,$ then $z\in J_1=R.$ Thus $ C_{\Delta}(R)=R.$

\vu

Now we prove that $R^\af=C(\Delta).$ We have $R^\af\subseteq C(\Delta)=C_{\Delta}(\Delta)\subseteq C_{\Delta}(R)=R.$ Let $r\in C(\Delta).$ Then for $x\in J_g$ we have $rx=xr=\af_{g}(r1_{g\m})x,$ which implies  $r1_gx=1_g(rx)=\af_{g}(r1_{g\m})x.$ Notice that $D_g\cong {\rm End}_{D_g}(J_g),$    implies that $J_g$ is a faithful $D_g$-module.  Then  $r1_g=\af_{g}(r1_{g\m}),$ for all $g\in G,$ and therefore $r\in R^\af.$

\vu

Finally we show that $\Delta/R^\af$ is separable.  By \cite[Theorem 2.7.1]{DI}   it is enough to prove that $\Delta_\mathfrak{m}/(R^\af)_\mathfrak{m}$ is separable, for every maximal ideal $\mathfrak{m}$ of $R^\af.$ But if $R^\af$ is local, then $R$ is semi-local being a f.g.p. over  $R^\af.$ It follows that    ${\bf Pic}(R)=0,$ and hence  ${\bf Pic}(D_g)=0,$ for all $g\in G$ as $D_g$ is a direct factor of $R.$ Thanks to the third item of
Lemma~\ref{gaction} $f$ is locally trivial, i.e.   $f(g)_\mfp \cong (D_g)_\mfp, $ for any $\mfp \in {\rm Spec}(R),$
  and it follows  from  Lemma \ref{iguald} that $f(g) \in {\bf Pic}(D_g)=0,$ so that  $f(g)\cong D_g$ as  $D_g$-modules as well as $R$-modules. Then, Example \ref{cpp}  implies that $\Delta_\mathfrak{m}$ is a crossed product by a  twisted partial action, and consequently, $\Delta_\mathfrak{m}$ is  separable over $R^\af$ thanks to \cite[Proposition 3.2]{PS}.  \cua

\vu

\begin{remark}\label{Phi_0} It is readily seen that the ring isomorphism   $\alpha _g : D_{g\m} \to D_g$ gives
 an $R$-$R$-bimodule isomorophism $$ _g(D_{g\m})_I \cong   _I(D_g)_{g^{-1}},$$ where
$   _I(D_g)_{g^{-1}} = D_g,$ as sets, the left $R$-action on $   _I(D_g)_{g^{-1}}$ is $ (r,d) \mapsto r d $ and the right $R$-action is defined by $(d, r) \mapsto d \alpha _g(r),$ $a\in  D_g, r\in R.$ Then obviously
$\Phi_0 (g ) = [_g(D_{g\m})_I] = [ _I(D_g)_{g^{-1}}] \in {\bf PicS}_{R^\af}(R).$
\end{remark}

 \begin{prop}\label{AisomorfoDelta}
		Let $A$ be an Azumaya $R^\af$-algebra containing $R$ as a maximal commutative subalgebra. Take  $f \in Z^1(G,\af^*,{\bf PicS}(R))$ such that  $\varphi_5([A])=cls(f)$.  Write $\Phi_f(g)=[J_g],$ then  there is a factor set $\mathcal{F}=\{f_{g,h}\colon J_g\ot J_h\to D_g\ot J_{gh}\}$ for $\Phi_f$  and  $\Delta=\Delta(\mathfrak{F},\af,R,\Phi_f, G)=\bigoplus_{g \in G}J_g$  is a partial generalized crossed product.
	\end{prop}
 \p Write $f(g)=[C_g]$ for all $g\in G,$ then  $J_g\simeq C_g\otimes  _I(D_g)_{g^{-1}} .$  By Remark \ref{mgtocg}  one may  assume that
$$C_g=\{a \in 1_{g\m} A \mid \alpha_{g\m}(r1_g)a=ar, \ \mbox{for any}\ r \in R\},$$
with left and right $R$-module actions given by $r\bullet a=\alpha_{g\m}(r1_g)a=ar=a\bullet r,$ for $r\in R$ and $a\in C_g.$
Let $J_g'=C_g$ as sets and  suppose that  the left and right $R$-module actions  on $J_g'$ are  given by

	\begin{equation}
	r\star a=\alpha_{g\m}(r1_g)a \ \ \mbox{and} \ \ a\star r=a\alpha_g(r1_{g\m}), \label{acoesdeJg}
	\end{equation}
	for $r\in R$ and $a\in J'_g.$
	
Consider the $R$-$R$-bimodule  map
	$$\begin{matrix}
	\lambda: & J_g' & \longrightarrow & C_g\otimes  _I(D_g)_{g^{-1}}  \\
	 & u_g& \longmapsto & u_g \otimes 1_g
	\end{matrix},$$

and let
$$\begin{matrix}
\zeta: & C_g\ot  _I(D_g)_{g^{-1}} & \longrightarrow & J_g',\\
 & c_g \ot d_g & \longmapsto & c_gd_g.
\end{matrix}$$ Then it is not difficult to check that $\zeta=\lambda\m,$ and $\lambda$  is an  isomorphism, so $J'_g\simeq J_g$ when considered with the  left and right $R$-module actions given by   (\ref{acoesdeJg}).

On the other hand, by Proposition \ref{pr1}   the map $\Phi_f$ is a partial representation, and we shall construct a factor set for it.

Let $u_g\in J_g$, $u_h \in J_h$ and $r \in R$, then
\begin{eqnarray*}
\af_{(gh)\m}(r1_{gh})u_hh_g & = & \af_{(gh)\m}(r1_{gh})1_{h\m}u_hh_g\\
& = & \af_{h\m}(\af_{g\m}(r1_g)1_h)u_hu_g\\
& = & u_h\af_{g\m}(r1_g)u_g\\
& = & u_hu_gr.
\end{eqnarray*}
In particular,
$$1_{(gh)\m}u_hu_g=\af_{(gh)\m}(1_g1_{gh})u_hu_g=u_hu_g1_g=u_hu_g.$$
Thus $u_hu_g \in J_{gh}$ and the map
$$\begin{matrix}
f_{g,h}: & J_g\otimes J_h & \longrightarrow & D_g\ot J_{gh}\\
& u_g\ot u_h & \longmapsto & 1_g\ot u_hu_g
\end{matrix},$$
is well defined, and it is not difficult to show that  $f_{g,h}$ is $R$-bilinear. Now we check that $f_{g,h}$ is an  isomorphism. Localizing by an ideal in $Spec(R^\af)$ we may assume that  $R^\af$  is a local ring.
  This implies that $R$ is semilocal and  $\Pic(R)=\{1\}$. Since $D_g$ is a direct summand of $R$,  we obtain $\Pic(D_g)=\{1\}$. But $[C_g]\in \Pic(D_g)$,  and there is a $R$-$R$-bimodule isomorphism $\gamma_g:D_g\longrightarrow C_g.$    Set $\omega_g=\gamma_g(1_g) \in C_g$,  then given  $c_g \in C_g$, there exists $d \in D_g$ such that
$$c_g=\gamma_g(d)=d\cdot \gamma_d(1_g)=\af_{g\m}(d)\omega_g=\omega_gd.$$
  This implies  that   $C_g=D_g\cdot \omega_g=\omega_gD_g$ with $\af_{g\m}(r1_{g\m})\omega_g=\omega_gr$, for each $r \in R$.

Consider the  isomorphism $\gamma_g'\colon _I(D_g)_{g^{-1}}  \to J_g$ given  by  the composition
$$ _I(D_g)_{g^{-1}}  \longrightarrow D_g\ot  _I(D_g)_{g^{-1}}  \longrightarrow
C_g\ot  _I(D_g)_{g^{-1}}  \longrightarrow  J_g.$$
 Then  $\gamma_g'(1_g)=\omega_g,$ and $J_g=D_{g\m}\omega_g=\omega_gD_g.$ From this we get that

$$\begin{matrix}
	f_{g,h}: & \om_gD_g\otimes \om_hD_h & \longrightarrow & D_g\ot \om_{gh}D_{gh}\\
	& \om_gd_g\ot \om_hd_h & \longmapsto & 1_g\ot \om_hd_h\om_gd_g
\end{matrix},$$
for all $g,h\in G.$
Now we construct the inverse of $f_{g,h}$.  By Lemma \ref{kan}, there is an  isomorphism $_g(1_{g\m}A)\simeq C_g\ot A$. Localizing by an ideal in $Spec(R^\af)$, we have
$_g(1_{g\m}A)\simeq \omega_gD_g\ot A\simeq \om_gA.$ Since $1_{g\m}\in 1_{g\m}A$,   there is  $a_g \in A$
 such that
\begin{equation}
\om_ga_g=1_{g\m}. \label{wgag}
\end{equation}
But
$$\omega_g(a_g\omega_g)=(\om_ga_g)\om_g=1_{g\m}\om_g=\om_g1_g.$$
Since  $\om_g$  is a free basis of  $_g(1_{g\m}A),$  we obtain
\begin{equation}
a_g\om_g=1_g. \label{agwg}
\end{equation}
For $r \in R$, one has that
\begin{eqnarray*}
\om_g\af_g(r1_{g\m})a_g & = & r\om_ga_g = r1_{g\m}\\
& = & 1_{g\m}r = \om_ga_gr.
\end{eqnarray*}
Then,
\begin{equation}
\af_g(r1_{g\m})a_g=a_gr, \ \ \mbox{for all }\ r \in R. \label{algag}
\end{equation}
But
\begin{eqnarray*}
1_ga_ga_h1_{h\m}\om_{gh}r & = & 1_ga_ga_h1_{h\m}\af_{h\m g\m}(r1_{gh})\om_{gh}\\
& = & 1_ga_ga_h1_{h\m}\af_{h\m}(\af_{g\m}(r1_{g})1_{h\m})\om_{gh}\\
&\stackrel{(\ref{algag})}{=} & 1_ga_g\af_{g\m}(r1_g)a_h1_{h\m}\om_{gh}\\
&\stackrel{(\ref{algag})}{=} & r1_ga_ga_h1_{h\m}\om_{gh}.
\end{eqnarray*}
Then, $1_ga_ga_h1_{h\m}\om_{gh} \in C_A(R)=R$. Hence,
$1_ga_ga_h1_{h\m}\om_{gh}=1_g1_{gh}a_ga_h\om_{gh} \in D_gD_{gh}$. Finally, we see that
\begin{eqnarray*}
\om_h\om_g(1_ga_ga_h1_{h\m}\om_{gh}) & = & \om_ga_ga_h1_{h\m}\om_{gh}\\
& \stackrel{(\ref{wgag})}{=} & \om_h1_{g\m}a_h1_{h\m}\om_{gh}\\
& = & \af_{h\m}(1_{g\m}1_h)\om_ha_h1_{h\m}\om_{gh}\\
& \stackrel{(\ref{wgag})}{=} & 1_{h\m g\m}1_{h\m}\om_{gh}\\
& = & 1_{h\m}\om_{gh}.
\end{eqnarray*}
We consider the well defined map
	$$\begin{array}{c c c l}
l_{g,h}: & D_g\ot \om_{gh}D_{gh} & \longrightarrow & \om_gD_g\ot \om_hD_h\\
&d_g\ot \om_{gh}d_{gh} & \longmapsto & \om_g(d_ga_ga_h1_{h\m}\om_{gh}d_{gh})\ot \om_h.
\end{array}$$

Moreover,
\begin{eqnarray*}
l_{g,h}\circ f_{g,h}(\om_gd_g\ot \om_hd_h) & = & l_{g,h}(1_g\ot \om_hd_h\om_gd_g)\\
& = &  l_{g,h}(1_g\ot \om_h\om_g\af_g(d_h1_{g\m})d_g)\\
& = & \om_g(1_ga_ga_h1_{h\m}\om_h\om_g\af_g(d_h1_{g\m})d_g)\ot \om_h\\
& = & \om_ga_ga_h\om_h\om_gd_g\af_g(d_h1_{g\m})\ot \om_h\\
& \stackrel{(\ref{agwg}),(\ref{wgag})}{=} & 1_{g\m}1_h\om_gd_g\af_g(d_h1_{g\m})\ot \om_h\\
& = & \af_{g\m}(1_{gh}1_g)\om_gd_g\af_g(d_h1_{g\m})\ot \om_h\\
& = & \om_g1_{gh}d_g\af_g(d_h1_{g\m})\ot \om_h\\
& = & \om_gd_g\af_g(1_h1_{g\m})\af_g(d_h1_{g\m})\ot \om_h\\
& = & \om_gd_g\af_g(d_h1_{g\m})\ot \om_h\\
& = & \om_gd_g\cdot d_h\ot \om_h\\
& = & \om_gd_g\ot d_h\cdot \om_h\\
& = & \om_gd_g\ot \af_{h\m}(d_h)\om_h\\
& = & \om_gd_g\ot \om_hd_h,
\end{eqnarray*}

and

\begin{eqnarray*}
f_{g,h}\circ l_{g,h}(d_g\ot \om_{gh}d_{gh}) & = & f_{g,h}(\om_g(d_ga_ga_h1_{h\m}\om_{gh}d_{gh})\ot \om_h)\\
& = & 1_g\ot \om_h\om_gd_ga_ga_h1_{h\m}\om_{gh}d_{gh}\\
& = & 1_g\ot \om_h\af_{g\m}(d_g)\om_ga_ga_h1_{h\m}\om_{gh}d_{gh}\\
& \stackrel{(\ref{wgag})}{=} & 1_g\ot \om_h\af_{g\m}(d_g)1_{g\m}a_h1_{h\m}\om_{gh}d_{gh}\\
& = & 1_g\ot \af_{h\m}(\af_{g\m}(d_g))\om_ha_h1_{h\m}\om_{gh}d_{gh}\\
& \stackrel{(\ref{wgag})}{=} & 1_g\ot \af_{h\m}(\af_{g\m}(d_g))1_{h\m}1_{h\m}\om_{gh}d_{gh}\\
& = & 1_g\ot \af_{h\m g\m}(d_g 1_{gh})\om_{gh}d_{gh}\\
& = & 1_g\ot d_g \cdot \om_{gh}d_{gh}\\
& = & d_g\ot \om_{gh}d_{gh}.
\end{eqnarray*}

\noindent and we conclude that $f_{g,h}$ is an isomorphism.
Thus to show that  $\mathcal{F}=\{f_{g,h}: J_g\ot J_g\longrightarrow D_g\ot J_{gh}, \ g,h \in G\}$ is a factor set,  it only remains to  prove that the diagram

$$\xymatrix{J_g\ot J_h\ot J_l\ar[rr]^{{\rm id}\ot f_{h,l}}\ar[dd]_{f_{g,h}\ot {\rm id}} & &J_g\ot D_h\ot J_{hl}\ar[rr]^{\varphi\ot {\rm id}} & & D_{gh}\ot J_g\ot J_{hl}\ar[d]^{{\rm id}\ot f_{g,hl}}\\
& & 	& &  D_{gh}\ot D_g\ot J_{ghl}\ar[d]^{m\ot {\rm id}}\\
D_g\ot J_{gh}\ot J_l\ar[rrrr]_{{\rm id}\ot f_{gh,l}} & &  & & D_g\ot D_{gh}\ot J_{ghl}},$$
where
$$\begin{matrix}
\varphi: & J_g\ot D_h & \longrightarrow & D_{gh}\ot J_g\\
& u_g \ot d_h & \longmapsto & \af_g(d_h1_{g\m})\ot u_g
\end{matrix}$$
and $$\begin{matrix}
m: & D_{gh}\ot D_{g}& \longrightarrow & D_g\ot D_{gh}\\
& d_{gh}\ot d_g & \longmapsto &d_g\ot d_{gh}
\end{matrix},$$
 is commutative.
 Indeed, for  $u_g \in J_g,u_h\in J_h, u_l \in J_l$, we have
\begin{eqnarray*}
u_g\ot u_h\ot u_l & \stackrel{f_{g,h}\ot {\rm id}}{\longmapsto} & 1_g\ot u_hu_g\ot u_l\\
& \stackrel{{\rm id}\ot f_{gh,l}}{\longmapsto}& 1_g\ot 1_{gh}\ot u_lu_hu_g.
\end{eqnarray*}
and,
\begin{eqnarray*}
u_g\ot u_h\ot u_l & \stackrel{{\rm id}\ot f_{h,l}}{\longmapsto} & u_g\ot 1_h\ot u_lu_h\\
& \stackrel{\varphi\ot {\rm id}}{\longmapsto}& 1_{gh}\ot u_g\ot u_lu_h\\
& \stackrel{D_{gh}\ot f_{g,hl}}{\longmapsto} & 1_{gh}\ot 1_g\ot u_lu_hu_g\\
& \stackrel{m\ot {\rm id}}{\longmapsto} & 1_g\ot 1_{gh}\ot u_ly_hu_g.
\end{eqnarray*}
then the diagram is commutative and $\Delta=\bigoplus_{g \in G}J_g$  is a partial generalized crossed product.
\cua

\begin{lema}\label{varphi5Delta}
Let  $f\in Z^1(G,\af^*,\Pics(R))$ be such that  $\Delta=\Delta(\mathcal{F}, \af, R,f\Phi_0,G)=\bigoplus_{g \in G}J_g$ is a partial generalized crossed product.  Then $\varphi_5([\Delta])=cls(f\m)$.
\end{lema}
\p
	Write$f(g)=[M_g]$, then $J_g\simeq M_g\ot  _I(D_g)_{g^{-1}} $. Since $M_g$ is a  $D_g$-module, we have that
	$J_g\ot  _I(D_{g\m})_{g}  \simeq M_g\ot  _I(D_g)_{g^{-1}}  \ot _I(D_{g\m})_g \simeq M_g\ot D_g\simeq M_g$. Hence,
	\begin{eqnarray}
	M_g\simeq J_g\ot   _I(D_{g\m})_g.  \label{MgisoJgDg\m}
	\end{eqnarray}
	By Proposition \ref{gpcpisa} we have that $[\Delta]\in B(R/R^\af)$ and the equality
	$$J_g=\{a \in \Delta: \af_g(r1_{g\m})a=ar, \ \mbox{for all}\ r\in R \},$$ follows by \eqref{carjh}.
	Denote $\varphi_5(\Delta)=cls(f')$, where $f'(g)=[C_g]$. By Lemma \ref{kan} we   obtain
	$$C_g=\{ a\in 1_{g\m}\Delta; \ \af_{g\m}(r1_g)a=ar, \ \mbox{for all } r\in R\}.$$
	Consider the well defined $R$-bimodule map
	$$\begin{array}{c c c l}
	\varphi: &   _I(D_g)_{g^{-1}}  \ot J_{g\m} & \longrightarrow & C_g\\
	& d\ot x_{g\m} & \longmapsto & \af_{g\m}(d)x_{g\m}.
	\end{array}$$
Then
	$$\begin{matrix}
	\xi: & C_g & \longrightarrow &  _I(D_g)_{g^{-1}}   \ot J_{g\m}\\
	& c_g & \longmapsto & 1_g\ot c_g
	\end{matrix}$$
is the  inverse of $\varphi$. Indeed, for  $c_g \in C_g$ we have
$$(\varphi\circ \xi)(c_g)=\varphi(1_g\ot c_g)=\af_{g\m}(1_g)c_g=1_{g\m}c_g=c_g.$$
Moreover, for  $d \in D_g$ and $x_{g\m}\in J_{g\m}$, we  see that
\begin{eqnarray*}
(\xi\circ \varphi)(d\ot x_{g\m}) & = & \xi(\af_{g\m}(d)x_{g\m})\\
& = & 1_g\ot x_{g\m}\\
& = & 1_g\cdot \af_{g\m}(d)\ot x_{g\m}\\
& = & d\ot x_{g\m},
\end{eqnarray*}
and we get that $\varphi=\xi\m,$ as  there is an isomorphism
\begin{equation}
C_g\simeq  _I(D_g)_{g^{-1}}  \ot J_{g\m}. \label{CgisoJgm}
\end{equation}
Since $f \in Z^1(G,\af^*,\Pics(R))$, then
$$\af_g^*(f(h)[D_{g\m}])f(gh)\m f(g)=[D_g][D_{gh}], \ \ \forall \ g,h \in G.$$
In particular, taking $h=g\m$ and using the fact   $f(1)=[R]$, we  obtain
\begin{equation}
f(g)\m=\af_g^*(f(g\m)[D_{g\m}]), \ \forall \ g\in G. \label{fg\m}
\end{equation}
Then
\begin{eqnarray*}
C_g& \simeq  & C_g\ot D_g \stackrel{(\ref{CgisoJgm})}{\simeq}   _I(D_g)_{g^{-1}}   \ot J_{g\m}\ot D_g\\
& \simeq &  _I(D_g)_{g^{-1}}  \ot J_{g\m} \ot  _I(D_g)_{g^{-1}}   \ot  _I(D_{g\m})_g\\
& \stackrel{(\ref{MgisoJgDg\m})}{\simeq} &  _I(D_g)_{g^{-1}}  \ot M_{g\m}\ot  _I (D_{g\m})_g\\
& \simeq & \af_g^*(f(g\m)[D_{g\m}])\\
& \stackrel{(\ref{fg\m})}{\simeq } & f(g)\m.
\end{eqnarray*}
Since $\varphi_5([\Delta])=cls(f')$, where $f'(g)=[C_g]=[f(g)\m]$, for all $g \in G$. Then $f'=f\m$ and $\varphi_5(\Delta)=cls(f\m)$.

\cua

\begin{teo} The sequence $B(R/R^\af)\stackrel{\varphi_5}{\to} H^1(G,\af^*, {\bf PicS}(R))\stackrel{\varphi_6}{\to}  H^3 (G,\af, R) $ is exact.
\end{teo}	
\p Let  $f \in Z^1(G,\af^*,\Pics(R))$  be such that  $cls(f)\in Im(\varphi_5)$, then there is  $[A]\in B(R/R^\af)$ with $\varphi_5([A])=cls(f)$. Write $f(g)=[C_g]$, it follows from Proposition \ref{AisomorfoDelta}  that $\Delta=\bigoplus_{g \in G}J_g$, where  $J_g=C_g\ot  _I(D_g)_{g^{-1}} $  is a partial generalized crossed product. Then, $\varphi_6(cls(f))$ is trivial in $H^3(G,\af,R)$ and consequently $cls(f)\in \ker(\varphi_6)$, and we obtain that $Im(\varphi_5)\subseteq \ker(\varphi_6)$.
	
	On the other hand, if $f\in Z^1(G,\af^*,\Pics(R))$ satisfies $cls(f)\in \ker(\varphi_6)$, then $\Delta=\bigoplus_{g \in G}J_g$, where  $J_g=f(g)\ot  _I(D_g)_{g^{-1}} $,  is a partial generalized crossed product and $[\Delta]\in B(R/R^\af)$. By  Lema \ref{varphi5Delta}, we have that  $\varphi_5(\Delta)=cls(f\m)$. Then $cls(f\m) \in Im(\varphi_5)$. Since $Im(\varphi_5)$ is a group we  conclude that  $cls(f)\in Im(\varphi_5),$  proving that  the sequence is exact.\cua

\vt

Summarizing, we have obtained  the following seven-term exact sequence:

\vd

$0\to H^1(G,\af, R)\stackrel{\f_1}{\to} {\bf Pic}(R^\af)\stackrel{\f_2}{\to }{\bf PicS}(R)^{\af^*}\cap {\bf Pic}(R) \stackrel {\f_3}{\to }H^2(G,\af, R) \stackrel{\f_4}{\to} B(R/R^\af)\stackrel{\f_5}{\to} H^1(G,\af^*,{\bf PicS}(R))\stackrel{\f_6}{\to}  H^3 (G,\af, R).$

\vd

\begin{cor} {\bf Hilbert's $90^{\rm th}$ Theorem for partial actions.} Let $R\supseteq R^\af$ be a partial Galois extension. If ${\bf Pic}(R^\af)=0,$ then $H^1(G,\af, R)=0.$ \cua
\end{cor}

\begin{cor}{\bf Crossed Product Theorem for partial actions.} Let $R$ be a commutative ring such that ${\bf Pic}(R)$ is trivial. If $\af$ is a unital partial action of a finite group $G$ on R such that $R/R^\af$ is an $\af$-partial Galois extension, then there is a group isomorphism $H^2(G,\af, R) \cong B(R/R^\af)$ given by ${\rm cls}(\om)\mapsto [R\star_{\af,\om} G].$\end{cor}
\p Since ${\bf Pic}(R)$ is trivial  we conclude that ${\bf Pic}(D_g)$ is trivial for all $g\in G.$ Then by the  third item of
Lemma~\ref{gaction}  we obtain that $f(g)=[D_g]$ for all $f\in Z^1(G,\af^*, {\bf PicS}(R))$ which implies  that $H^1(G,\af^*,{\bf PicS}(R))$ is trivial and $\f_4$ is an isomorphism. \cua
\begin{cor}\cite {CHR} Let G be a finite group acting globally on the commutative ring R and suppose that $R/R^G$ is a Galois extension. Then there is a seven term exact sequence
$0\to H^1(G, \U(R))\stackrel{\f_1}{\to} {\bf Pic}(R^G)\stackrel{\f_2}{\to } {\bf Pic}(R)^G \stackrel {\f_3}{\to }
H^2(G, \U(R)  ) \stackrel{\f_4}{\to} B(R/R^\af)\stackrel{\f_5}{\to} H^1(G,{\bf Pic}(R))\stackrel{\f_6}{\to} H^3(G, \U(R)).$
\end{cor}

\p  When the action $\af$ is global we have $R^\af=R^G$ and $D_g=R$ for any $g\in G.$ The latter implies   $H^i(G,\af, R)=H^i(G, \U(R)),$ for $1\le i\le 3.$
Moreover ${\bf PicS}(R)^{\af^*}={\bf PicS}(R)^{G}=\{[E]\in {\bf PicS}(R)\mid E_g\cong E\,\text {for any}\, g\in G \}$ and then ${\bf PicS}(R)^{\af^*}\cap {\bf Pic}(R)={\bf Pic}(R)^G.$  Finally,
\begin{align*}{\rm cls} (f) \in   H^1 (G,\af^*,{\bf PicS}(R)) & \Leftrightarrow
f\in Z^1(G,\af^*,{\bf PicS}(R)) \,\, \&\,\,   f(g) \in {\bf Pic}(D_g), \, \forall  g\in G.\\
& \Leftrightarrow  f\in Z^1(G,\af^*,{\bf PicS}(R)) \,\, \&\,\, f(g) \in {\bf Pic}(R), \, \forall  g\in G.\\
& \Leftrightarrow  f\in Z^1 (G,{\bf Pic}(R)),
\end{align*}
so that   $ H^1 (G,\af^*,{\bf PicS}(R)) =H^1 (G,{\bf Pic}(R)),$
and we obtain the sequence as given in \cite {CHR}.   \cua

\begin{remark} The fact that  $ \varphi_6 (cls(f)) = cls(\om ) \in  H^3(G,\af,R)$  for every $1$-cocycle $f\in Z^1(G,\af^*, {\bf PicS}(R))$ was proved in \cite[Section 6.1]{DPP1} by showing that $\om  _\mfp $ is a $3$-coboundary  for any prime $\mfp\in \rm{Spec}(R^{\alpha}).$ Thus it is reasonable  to say that a partial cocycle  $\om\in Z^n(G,\af,R)$ is a local coboundary if for any $\mfp\in \rm{Spec}(R^{\alpha})$ there exists
 $\rho\in C^{n-1}(G,\af  _\mfp , R_ \mfp )$ such that  $\om(g_1, \ldots ,g_n)  _\mfp=(\de^{n-1}\rho)(g_1 ,\ldots , g_n),$ where  $\af  _\mfp $ is the partial action of $G$ on $R _\mfp $ which corresponds to $\alpha .$ Then we may call the cohomology class of   $\om\in Z^n(G,\af,R)$   locally trivial if  $\om $ is a local coboundary, and we denote by $H^n_{\rm lt}(G,\af,R)$ the subgroup of $H^n (G,\af,R)$ which consists of the locally trivial cohomology classes. With this notation we have that  ${\rm Im} \, \varphi_6 \subseteq  H^3_{\rm lt}(G,\af,R).$ \end{remark}

\end{document}